\documentclass[a4paper, oneside, 11pt]{amsart}
\usepackage[utf8]{inputenc}

\usepackage{amsfonts}
\usepackage{amsmath}
\usepackage{amssymb}
\usepackage{amsthm}
\usepackage[a4paper]{geometry}
\usepackage{mathrsfs} 
\usepackage[dvipsnames]{xcolor}

\usepackage[pagebackref=true,colorlinks=true, linkcolor=MidnightBlue, citecolor=MidnightBlue]{hyperref}
\renewcommand*\backref[1]{\ifx#1\relax \else (Cited on #1) \fi}

\usepackage{appendix}
\usepackage{enumerate}
\usepackage{tikz-cd} 

\theoremstyle{plain}
\newtheorem{theorem}{Theorem}
\newtheorem{proposition}[theorem]{Proposition}
\newtheorem{lemma}[theorem]{Lemma}
\newtheorem{corollary}[theorem]{Corollary}
\newtheorem{definition}[theorem]{Definition}
\newtheorem{remark}[theorem]{Remark}

\theoremstyle{definition}

\numberwithin{theorem}{section}
\numberwithin{equation}{section} 

\newcommand{\norm}[1]{\left \lVert  #1 \right \rVert}
\newcommand{\abs}[1]{\left\lvert #1 \right\rvert}

\newcommand{\vertiii}[1]{{\left\vert\kern-0.25ex\left\vert\kern-0.25ex\left\vert #1 
    \right\vert\kern-0.25ex\right\vert\kern-0.25ex\right\vert}}

\newcommand{\Z}{\mathbb{Z}}

\newcommand{\R}{\mathbb{R}}
\newcommand{\N}{\mathbb{N}}

\newcommand{\calF}{\mathcal{F}}

\newcommand{\calH}{\mathcal{H}}

\newcommand{\calM}{\mathcal{M}}

\title[Time-periodic behaviour in IPS]{Time-periodic behaviour in one- and two-dimensional interacting particle systems}
\author{Benedikt Jahnel}
\address{Institut f\"ur Mathematische Stochastik, Technische Universit\"at Braunschweig, Universit\"atsplatz 2,
38106 Braunschweig, Germany \& Weierstrass Institute for Applied Analysis and Stochastics\\
Mohrenstraße 39\\
10117 Berlin\\
Germany}
\email{benedikt.jahnel@tu-braunschweig.de}
\author{Jonas Köppl}
\address{Weierstrass Institute for Applied Analysis and Stochastics\\
Mohrenstraße 39\\
10117 Berlin\\
Germany}
\email{koeppl@wias-berlin.de}
\date{\today}

\keywords{Interacting particle systems, Gibbs measures, periodic behaviour, attractor properties, synchronisation, non-ergodicity, non-equilibrium stationary state}
\subjclass{Primary 82C22; Secondary 60K35} 

\begin{document}

\maketitle

\begin{abstract}
    We provide a class of examples of interacting particle systems on $\Z^d$, for $d\in\{1,2\}$, that admit a unique translation-invariant stationary measure, which is not the long-time limit of all translation-invariant starting measures, due to the existence of time-periodic orbits in the associated measure-valued dynamics. 
    This is the first such example and shows that even in low dimensions, not every limit point of the measure-valued dynamics needs to be a time-stationary measure. 
\end{abstract}

\section{Introduction and motivation}\label{sec:introduction}
We consider interacting particle systems on $\Z^d$ as in \cite{liggett_interacting_2005}, i.e., we study Markov processes on a state space of the form $\Omega = \Omega_0^{\Z^d}$ for some finite set $\Omega_0$ called the local state space. These processes are usually specified via their generator $\mathscr{L}$, given in terms of time-homogeneous transition rates, and its associated Markovian semigroup $(S(t))_{t \geq 0}$ acting on the space of continuous functions $C(\Omega)$. More precisely, the generator is usually of the form 
\begin{align*}
    \mathscr{L}f(\eta) = \sum_{x\in \Z^d}\sum_{\xi_x \in \Omega_0}c_x(\eta,\xi_x)\left[f(\xi_x\eta_{x^c})-f(\eta)\right], 
\end{align*}
where $c_x(\eta,\xi_x)\ge 0$ should be interpreted as the infinitesimal rate at which the particle at site $x \in \Z^d$ switches from the state $\eta_x$ to $\xi_x$, given that the rest of the system is currently in state $\eta_{x^c}$. 
In this setting, one can ask the following question. 
\begin{enumerate}[\bfseries (Q1)]
    \item Let $\mu_0$ be a probability measure on $\Omega$ and $\mu_t$ the distribution of the process with generator $\mathscr{L}$ at time $t \geq 0$. Under which conditions on the generator $\mathscr{L}$ are all weak limit points $\mu^*$ of $(\mu_t)_{t \geq 0}$ stationary measures for the process?
\end{enumerate} 
To exclude trivial non-random counterexamples we focus our attention on interacting particle systems which are \textit{non-degenerate} in the sense that they are truly random and (locally) irreducible:
\begin{enumerate}[\bfseries{(ND)}]
    \item For every finite volume $\Delta \Subset \Z^d$ and every pair of configurations $\eta, \xi \in \Omega$, that agree on $\Delta^c$, there exists a finite sequence $\eta^{(0)}, \dots, \eta^{(n)}$ of configurations such that $\eta^{(0)}=\eta$, $\eta^{(n)} = \xi$ and the transition rate from $\eta^{(i)}$ to $\eta^{(i+1)}$ is positive for all $i \in \{0,\dots, n-1\}$.
\end{enumerate}
This implies in particular, that for every time $T>0$, the transition probability $p_T(A,B)$ for any two local events $A,B \subset \Omega$ is strictly positive, see \cite[Proposition 6.1]{jahnel_long-time_2023} for a proof of this in a similar situation. 

While a classical result due to Mountford and Ramirez--Varadhan, see \cite{mountford_coupling_1995} and \cite{ramirez_relative_1996}, from the 90s shows that, for one-dimensional interacting particle systems with finite-range interactions, all limit points of the associated measure-valued dynamics are stationary measures, essentially nothing was known about what can or cannot happen in dimensions $d\geq 2$.

Perhaps the simplest way in which a limit point $\mu^*$ of $(\mu_t)_{t \geq 0}$ can fail to be a stationary measure is if the system admits a \textit{time-periodic orbit}, i.e., there exists a probability measure $\mu_0$ that has the property that, if the system at time zero is distributed according to this law, then it returns to this law after a finite time $T>0$, but the system has a different distribution at all intermediate times $0<t<T$.
While it is not so hard to show that this behaviour cannot happen in continuous-time Markov chains on finite state spaces, the situation for interacting particle systems on the infinite lattice $\Z^d$ is much more delicate and very little is known rigorously about interacting particle systems with periodic laws.

The first example of a non-degenerate interacting particle system which exhibits time-periodic behaviour was constructed in \cite{jahnel_class_2014}, but due to regularity problems, this construction, which was inspired by conjectures in \cite{maes_rotating_2011}, only works in dimensions $d\geq 3$ and left the question of the possibility of time-periodic behaviour for two dimensions in general, and systems with long-range interactions in one spatial dimension widely open. 

More recently, it was shown in \cite{jahnel_long-time_2023} that in dimensions $d=1,2$ reversible interacting particle systems with short-range interactions cannot exhibit time-periodic behaviour. So to construct interacting particle systems with time-periodic behaviour in dimensions one and two, it is necessary to violate one of these two assumptions.

In this article, we provide a new proof strategy that allows to extend the construction in \cite{jahnel_class_2014} to one and two dimension.
In particular, we show that even low-dimensional systems can exhibit time-periodic behaviour and thereby also establish that the result of Mountford and Ramirez--Varadhan does not extend to interacting particle systems with arbitrary interaction range. More precisely, we prove the following. 

\begin{theorem}\label{theorem:main-result}
    There exists a non-degenerate interacting particle system on $\Z^d$ for $d=1,2$, with finite local state space, translation-invariant and time-homogeneous transition rates and Markov semigroup $(S(t))_{t\ge 0}$ that satisfies the following properties. 
    \begin{enumerate}[(a)]
        \item It admits a unique translation-invariant time-stationary measure. 
        \item There exists a non-trivial family of measures $\{\nu_t \colon t \in [0,2\pi]\}$ that are not time-stationary w.r.t.~the dynamics but such that 
        \begin{align*}
            \forall s,t \geq 0: \quad \nu_tS(s) = \nu_{t+s \ \text{mod } 2\pi}.
        \end{align*}
    \end{enumerate}
\end{theorem}

This improves and extends the construction from \cite{jahnel_class_2014}, which only works in dimensions $d \geq 3$, since it, on the one hand, needs continuous symmetry breaking but, on the other hand, also a sufficiently fast decay of the interactions between particles. To be precise, if one checks their arguments carefully it requires that the oscillations of the rates decay faster than $\abs{x}^{-2d}$, as $x \to \infty$. In dimensions $d \geq 3$ this is not a problem, because there, even continuous systems with nearest-neighbour interactions exhibit first-order phase transitions. But this is not true in dimensions one and two because there, only systems with long-range interactions that decay slower than $\abs{x}^{-2d}$, as $x \to \infty$, can exhibit continuous symmetry breaking. 
Despite the additional difficulties of long-range dependencies, our new proof strategy is based on more elementary ideas and can also be used for short-range systems in $d\geq3$. Therefore, we believe that it also makes the main ideas behind the original construction more transparent without obfuscating them with the technical details of their implementation.

For the moment, we only know that the time-periodic orbit exists but we do not know if it is stable in the sense of being (locally) attractive. We do know, however, that it is somewhat stable if we perturb the dynamics a little bit. See Remark \ref{remark:stability-under-reversible-perturbations} for a precise statement of this stability property.  

Let us also note that if one considers infinite local state spaces, e.g., the unit circle, and Langevin-type dynamics, one can construct an easier example of a system of interacting diffusions exhibiting non-trivial time-periodic behaviour, see \cite[Section 4]{maes_rotating_2011}. However, to our knowledge, there is to date no counterpart of the no-go results from \cite{mountford_coupling_1995}, \cite{ramirez_relative_1996}, and \cite{jahnel_long-time_2023} in this continuous setting. Testing the limitations of these no-go theorems in their natural setting requires a more involved construction, see Section \ref{section:brief-overview}, which is however still inspired by the continuous situation.

\subsection{Organisation of the manuscript}
The rest of this article is organised as follows. We first give a brief non-technical description of the construction in Section \ref{section:brief-overview}. 
Afterwards, we provide a completely rigorous roadmap for the construction that is necessary to prove our main result in  Section \ref{section:technical-overview}. The proofs can then be found in Section \ref{section:proofs}. We conclude with a brief outlook and some questions for further research in Section \ref{section:outlook}.

\section{Construction overview}\label{section:brief-overview}
Let us quickly sketch the main idea behind the construction in a non-technical way to avoid getting lost in the technical details that are necessary for a rigorous proof of Theorem \ref{theorem:main-result}. 
Briefly summarised, we proceed as follows. We consider a system with continuous local state space $\Omega_0 = \mathbb{S}^1$, that exhibits continuous symmetry breaking, i.e., non-uniqueness of the Gibbs measure for sufficiently low temperatures. Since we are working in one and two dimensions, we cannot use the classical nearest-neighbour XY model, because it does not undergo a first-order phase transition in $d\leq 2$, see \cite[Chapter 9]{friedli_statistical_2017}. Instead, we consider a \textit{long-range} version of the XY model, which does exhibit a phase transition in low dimensions, if one makes the interaction strength decay sufficiently slow, see \cite{kunz_first_1976}. This is the content of Section \ref{section:gibbs-measures}. 
\medskip

In the system with continuous spins, we can now apply a deterministic rotation, i.e., for every $t \in \R$ let $R_t : \Omega_0 \to \Omega_0$ be the rotation by the angle $t$ in the clockwise direction. In particular, if we apply this rotation to every site, then this induces a rotation of the extremal Gibbs measures associated to the long-range interaction potential. Of course, this deterministic dynamics has periodic orbits with period given by $2\pi$.\medskip

From a stochastic point of view this is Markovian but very degenerate, because there is absolutely no randomness involved yet. But we can now discretise the local state space $\mathbb{S}^1$ into $q$ arcs of length $2\pi/q$ for $q \in \N$ sufficiently large. After this discretisation, the deterministic rotation corresponds – in some sense to be made precise – to a stochastic dynamic with finite local state space $\Omega'_0 = \{1,\dots, q\}$. A more detailed description on how this discretisation works and how one recovers a specification for the discretised Gibbs measures can be found in Section \ref{section:discretising-gibbs-measures}. 

The main task is now to identify a particular time-periodic measure-valued trajectory of the deterministic rotation with a trajectory of the measure-valued discretised system and to show that this trajectory is also the time-evolution of an appropriately chosen interacting particle system. By being very careful in the construction, one can then deduce that this interacting particle system admits a unique equilibrium measure and fails to be ergodic due to the existence of this time-periodic orbit. This is the most subtle part of the proof and a more detailed description of how this is done can be found in Section \ref{section:interacting-particle-systems}. 
\medskip

Making this last part work in the presence of long-range interactions is the main technical contribution of this article, because we need to deviate quite a bit from the original proof in \cite{jahnel_class_2014}, which dealt with the technically easier case of short-range interactions in $d\geq 3$. In particular, we provide a new way, which does not depend on the regularity of the transition rates (as long as they are well-defined), to identify measure-valued trajectories of the interacting particle system and the discretised rotation.
Since our new strategy is based on more elementary ideas and can also be used for short-range systems in $d\geq3$, we believe that this also makes the main ideas behind the original construction more transparent without obfuscating them with the technical details of their implementation.

\section{Strategy of proof}\label{section:technical-overview}
\subsection{Setting and notation}\label{section:setting-and-notation}
Let us consider the continuous local state space $\Omega_0 = \mathbb{S}^1 \simeq [0, 2\pi)$ and the configuration space $\Omega = \Omega_0^{\Z^d}$ for $d \in \{1,2\}$. Let $\lambda$ denote the Lebesgue measure on $\mathbb{S}^1$. Further, let $\calF$ denote the Borel sigma algebra based on the product topology on $\Omega$ and write $\calM_1(\Omega)$ for the space of probability measures on $(\Omega,\calF)$. For $\Lambda\subset \Z^d$ we write $\Omega_\Lambda=\Omega_0^\Lambda$,  $\Lambda^c=\Z^d\setminus \Lambda$ and denote by $\calF_\Lambda\subset\calF$ the sub-sigma algebra of events measurable w.r.t.~$\Lambda$. We will use the shorthand notation $\Lambda \Subset \Z^d$ to signify that $\Lambda$ is a finite subset of $\Z^d$. 
\subsection{Gibbs measures and the DLR formalism}\label{section:gibbs-measures}
For our construction we will need to construct special families of probability measures on $\Omega$, so called Gibbs measures, via the DLR formalism. Let us first recall the general definition of a specification. 
\begin{definition}
A \emph{specification} $\gamma = (\gamma_\Lambda)_{\Lambda \Subset \Z^d}$ is a family of probability kernels $\gamma_{\Lambda}$ from $\Omega$ to $\calM_1(\Omega)$ that additionally satisfies the following properties. 
\begin{enumerate}[i.]
    \item Each $\gamma_{\Lambda}$ is \emph{proper}, i.e., for all $B \in \calF_{\Lambda^c}$ it holds that 
    \begin{align*}
        \gamma_\Lambda(B | \cdot) = \mathbf{1}_B(\cdot). 
    \end{align*}
    \item The probability kernels are \emph{consistent} in the sense that if $\Delta \subset \Lambda \Subset \Z^d$, then for all $A \in \calF$
    \begin{align*}
        \gamma_{\Lambda}\gamma_{\Delta}(A|\cdot)
        = 
        \gamma_\Lambda(A|\cdot),
\end{align*}
where the concatenation of two probability kernels is defined as usual via 
\begin{align*}
    \gamma_\Lambda \gamma_\Delta(A \lvert \eta) = \int_\Omega \gamma_\Delta(A \lvert \omega) \gamma_\Lambda(d\omega \lvert \eta). 
\end{align*}
\end{enumerate}
An infinite-volume probability measure $\mu$ on $\Omega$ is called a \emph{Gibbs measure} for $\gamma$ if $\mu$ satisfies the so-called \emph{DLR equations}. That is, for all $\Lambda \Subset \Z^d$ and $B\in \calF$ we have the following disintegration property
\begin{align*}
    \mu(\gamma_\Lambda(B\lvert \cdot)) = \mu(B). 
\end{align*}
We will denote the set of all Gibbs measures for a specification $\gamma$ by $\mathscr{G}(\gamma)$. 
\end{definition}
For the existence and further properties of Gibbs measures with specification $\gamma$ one needs to impose some conditions on the specification $\gamma$. More precisely, we will need the following. 
\begin{definition}
    A specification $\gamma$ is called:
    \begin{enumerate}[(i)]
        \item \emph{Translation-invariant}, if for all $\Lambda \Subset \Z^d$ and $i \in \Z^d$ we have 
        \begin{align*}
            \gamma_{\Lambda +i}(\eta_{\Lambda+i}\lvert \eta_{(\Lambda+i)^c})=\gamma_{\Lambda}(\eta_{\Lambda}\lvert \eta_{\Lambda^c}),
        \end{align*}
        where $(\Lambda +i)$ denotes the lattice translate of $\Lambda$ by $i$. 
        \item \emph{Nonnull}, if for some $\delta >0$,
        \begin{align*}
            \inf_{\eta \in \Omega}\gamma_0(\eta_0 \lvert \eta_{0^c}) \geq \delta. 
        \end{align*}
        \item \emph{Quasilocal}, if for all $\Lambda \Subset \Z^d$,
        \begin{align*}
            \lim_{\Delta \uparrow \Z^d}\sup_{\eta, \xi \in \Omega}\abs{\gamma_\Lambda(\eta_\Lambda \lvert \eta_{\Delta \setminus \Lambda}\xi_{\Delta^c}) - \gamma_\Lambda(\eta_\Lambda \lvert \eta_{\Lambda^c})}=0. 
        \end{align*}
    \end{enumerate}
\end{definition}
One sufficient condition to ensure the existence of a Gibbs measure for a given specification $\gamma$ is quasilocality, see e.g., \cite{georgii_gibbs_2011} or \cite{friedli_statistical_2017}.
For example, specifications defined via a translation-invariant uniformly absolutely summable potential $\Phi =(\Phi_A)_{A \Subset \Z^d}$ are translation-invariant, nonnull and quasilocal. On an intuitive level, being nonnull means that the conditional probability of a fixed vertex $x$ being in state $\eta_x \in \Omega_0$ is bounded from below \textit{uniformly} in $\eta_x$ and the rest of the configuration $\eta_{x^c}$.

In the following, we will always work with a specification $\gamma$ that is given in terms of an interaction potential $\Phi = (\Phi_A)_{A \Subset \Z^d}$, i.e., it is given by 
\begin{align}\label{eqn:finite-volume-lebesgue-density}
    \gamma_{\Lambda,\beta}(B\lvert \eta) 
    = 
    \frac{1}{Z_{\Lambda,\beta}^\eta}\int_{\Omega_\Lambda}\mathbf{1}_B(\omega_\Lambda\eta_{\Lambda^c})\exp(-\beta\calH_\Lambda(\omega_\Lambda\eta_{\Lambda^c}))\lambda^{\otimes\Lambda}(d\omega_\Lambda) 
\end{align}
for a finite volume $\Lambda \Subset \Z^d$, a measurable set $B\in \calF$, a boundary condition $\eta \in \Omega$, and the Hamiltonian 
\begin{align*}
    \calH_\Lambda(\omega) = \sum_{A \cap \Lambda \neq \emptyset}\Phi_A(\omega).
\end{align*}
In that case, we will also write $\mathscr{G}(\beta\Phi)$ for the set of all Gibbs measures w.r.t.~this specification. Under some mild assumptions on $\Phi$ one can show that this set is non-empty and convex.  We will denote the set of its translation-invariant extremal points by $\text{ex}\mathscr{G}_\theta(\beta \Phi)$. 
\medskip 

For our construction, we work with a long-range version of the ferromagnetic Heisenberg model, i.e., the interaction potential $\Phi$ we consider is given by 
\begin{align*}
    \Phi_A(\omega) 
    =
    \begin{cases}
        \abs{x-y}^{-\alpha}\omega(x) \cdot \omega(y) \quad &\text{if } A = \{x,y\}, \\\
        0 &\text{otherwise,}
    \end{cases}
\end{align*}
where $\alpha \in (d,2d)$ and $u \cdot v$ denotes the scalar product of two vectors $u,v\in \R^2$. 
If we use the parametrisation $\mathbb{S}^1 \ni \omega \mapsto \varphi \in [0,2\pi)$, the scalar product $\omega_x \cdot \omega_y$ corresponds to $\cos(\varphi_x - \varphi_y)$. This is why this type of interaction is both known as scalar product interaction and cosine interaction in the literature. Note that this potential is absolutely summable and translation invariant and hence, the associated specification is translation invariant, nonnull and quasilocal. 

The reason why we are interested in this class of models is because they exhibit continuous symmetry breaking, even in dimensions one and two whereas finite-range models like the classical nearest-neighbour XY model do not have a first-order phase transition in dimensions $d \leq 2$, see \cite[Chapter 9]{friedli_statistical_2017}. The following result was first shown in \cite{kunz_first_1976} for dimension $d=2$.
\begin{theorem}[Symmetry breaking in the long-range Heisenberg model]\label{theorem:continous-symmetry-breaking}
    If $\beta >0$ is sufficiently large, then there exists a number $c(\beta)>0$ and a family $\{\mu_u: \ u \in \R^2, \abs{u} = c(\beta)\} \subset \emph{ex}\mathscr{G}_\theta(\beta\Phi)$ such that $\mu_u(\sigma_0) = u$. 
\end{theorem}
In words, at sufficiently low temperatures there is spontaneous magnetisation and thereby a breaking of the $SO(2)$-symmetry. For a proof see~\cite[Theorem 20.15 \& Example 20.21]{georgii_gibbs_2011}. 

From now on we will always assume that $\beta$ is sufficiently large such that we are in the regime of Theorem \ref{theorem:continous-symmetry-breaking} and drop it from our notation.

\subsection{Discretising Gibbs measures}\label{section:discretising-gibbs-measures}
The processes we will construct will have the discretised circle as local state space, i.e., the set $\mathbb{S}^1_q = \{1,\dots, q\}$ where $1 \leq l \leq q$ stands for the $l$-th arc in a equidistant partition of the continuous circle into $q$ arcs.
Denote the associated discretisation map (or coarse graining) by $T_q : \mathbb{S}^1 \to \mathbb{S}^1_q$. Denote the discrete configuration space by $\Omega'_q = (\mathbb{S}^1_q)^{\Z^d}$.
We will later need to choose the parameter $q$ of this coarse graining large enough so that the image measure of $\mu \in \mathscr{G}(\gamma)$ under this discretisation is again a Gibbs measure with respect to some discretised specification $\gamma'$. 

\subsubsection{Dobrushin's comparison theorem and its implications}
In order to define this specification, we will need to introduce some more machinery, in particular, we will first need a family of specifications that allows us to go back from the coarse-grained model to a model with continuous spins. For this we rely on some classical results from \cite[Chapter 8]{georgii_gibbs_2011}. For a specification $\gamma$ we define the \textit{Dobrushin interdependence matrix} $C$ by 
\begin{align*}
    C_{xy}(\gamma) = \sup_{\xi, \eta \in \Omega\colon \xi_{y^c} = \eta_{y^c}}\norm{\gamma_x(\cdot\lvert \xi) - \gamma_x(\cdot \lvert \eta)}_{\text{TV}},
\end{align*}
where $\norm{\cdot}_{\text{TV}}$ is the total variation distance. In order to define the discretised specification $\gamma'$ we will need to use the following result in an intermediate step. 

\begin{theorem}[Dobrushin's comparison theorem]\label{theorem:dobrushin-comparison}
    Let $\gamma$ and $\Tilde{\gamma}$ be two specifications and suppose that $\gamma$ satisfies \emph{Dobrushin's condition}, i.e., $\gamma$ is quasilocal and such that 
    \begin{align*}
        c(\gamma) = \sup_{x \in \Z^d}\sum_{y \in \Z^d}C_{xy}(\gamma) < 1. 
    \end{align*}
    For each $x \in \Z^d$ let $b_x:\Omega \to \R$ be a measurable function such that for all $\omega \in \Omega$
    \begin{align*}
        \norm{\gamma_x(\cdot \lvert \omega) - \Tilde{\gamma}_x(\cdot \lvert \omega)}_{\text{TV}} \leq b_x(\omega). 
    \end{align*}
    Then, for $\mu \in \mathscr{G}(\gamma)$ and $\Tilde{\mu} \in \mathscr{G}(\Tilde{\gamma})$ it holds that, 
    for all $f\in C(\Omega)$,
    \begin{align*}
        \abs{\mu(f) - \Tilde{\mu}(f)} \leq \sum_{x,y\in \Z^d}\delta_x(f) D_{x,y}(\gamma)\Tilde{\mu}(b_y), 
    \end{align*}
    where the $D_{xy}(\gamma)$ are entries of the matrix $D(\gamma)$ defined by
    \begin{align*}
        D(\gamma) = \sum_{n\geq 0}C^n(\gamma). 
    \end{align*}
Here, the oscillation of a function $f:\Omega \to \R$ at the site $x$ is defined as
    \begin{align*}
        \delta_x(f) = \sup_{\eta, \xi: \ \eta_{x^c} = \xi_{x^c}}\abs{f(\eta)-f(\xi)}.
    \end{align*}
\end{theorem}
By applying this in the situation where $\gamma = \Tilde{\gamma}$, one immediately gets uniqueness of the Gibbs measure for $\gamma$ if $c(\gamma)<1$. This corollary of the comparison theorem is often referred to as the \textit{Dobrushin uniqueness theorem}. 

We will mostly apply this to the family of specifications which is the content of the next lemma. For $\omega' \in \Omega'_q$ and $\Lambda \subset \Z^d$ we define 
\begin{align*}
    [\omega'_\Lambda] = \{\omega \in \Omega: \ T(\omega)\lvert_\Lambda = \omega'_\Lambda\}. 
\end{align*}
With this notation at hand, we can state the next step in our construction. 

\begin{lemma}\label{eqn:definition-constrained-specification}
    For each fixed discrete-spin configuration $\omega' \in \Omega'_q$, define a family of probability kernels on the continuous spin space by constraining the continuous spins to $[\omega']$ and putting, for each finite $\Lambda \Subset \Z^d$, and bounded measurable observable $f:\Omega \to \R$
    \begin{align*}
        \gamma^{\omega'}_\Lambda(f \lvert \omega_{\Lambda^c}) 
        =
        \frac{\gamma_\Lambda(f \mathbf{1}_{[\omega'_\Lambda]}\lvert \omega_{\Lambda^c})}{\gamma_\Lambda(\mathbf{1}_{[\omega'_\Lambda]}\lvert \omega_{\Lambda^c})}, \quad \omega \in \Omega. 
    \end{align*}
    Then, $\gamma^{\omega'} = (\gamma^{\omega'}_\Lambda)_{\Lambda \Subset \Z^d}$ defines a quasilocal specification on the continuous spin space $\Omega$. 
\end{lemma}
For a proof of this see~\cite[Lemma 2.2]{jahnel_class_2014}. By standard existence results for Gibbs measures, this implies that for every $\omega'$ the set $\mathscr{G}(\gamma^{\omega'})$ of Gibbs measures compatible with the specification $\gamma^{\omega'}$ is not empty. For our construction to work, we now need to make sure that for each $\omega'$ it consists of just one element, i.e., for every $\omega' \in \Omega'_q$ there is a unique Gibbs measure $\mu^{\omega'}$ for the specification $\gamma^{\omega'}$. For this we will make use of Theorem \ref{theorem:dobrushin-comparison}. To apply the estimate from Theorem \ref{theorem:dobrushin-comparison} we will consider \textit{uniform} Dobrushin coefficients for the whole family of specifications $(\gamma^{\omega'})_{\omega' \in \Omega'}$. These are defined by 
\begin{align*}
    \overline{C}_{xy} = \sup_{\omega' \in \Omega'_q}\sup_{\xi, \eta \in \Omega: \ \xi_{y^c} = \eta_{y^c} \text{ and } T(\eta) = T(\xi)=\omega'}\norm{\gamma_x(\cdot\lvert \xi) - \gamma_x(\cdot \lvert \eta)}_{\text{TV}}.
\end{align*}
Since we took the supremum over all $\omega'\in \Omega'_q$, the associated Dobrushin constant
\begin{align*}
    \overline{c} = \sup_{x \in \Z^d}\sum_{y \in \Z^d}\overline{C}_{xy}
\end{align*}
is uniform in $\omega'$. So if we can ensure that $\overline{c}<1$, then, by Dobrushin's uniqueness theorem, for every $\omega'$ there exists a unique Gibbs measure with respect to the specification $\gamma^{\omega'}$ and of course we can also use the estimate from Theorem \ref{theorem:dobrushin-comparison} uniformly. 

\subsubsection{Choosing the discretisation fine enough}
Intuitively one would expect that, in order to get $\overline{c}<1$, we will need to choose the discretisation fine enough, i.e., $q$ large enough. This intuition is made precise via the following result which first appeared in \cite[Theorem 2.1]{van_enter_discrete_2011}. To state the result we will first need to set up some more notation. For every pair of sites $x,y \in \Z^d$ define a distance for $\varphi_x, \xi_x \in \mathbb{S}^1$  via
\begin{align*}
    d_{xy}(\varphi_x, \xi_x)=
    &\hspace{-0.55cm}\sup_{\omega, \eta\colon\omega\lvert_{y^c}=\eta\lvert_{y^c}, T(\omega_y) = T(\eta_y)}
    \hspace{-0.55cm}\abs{\calH_x(\varphi_x \omega_{x^c})-\calH_x(\varphi_x\eta_{x^c})-(\calH_x(\xi_x\omega_{x^c})-\calH_x(\xi_x\eta_{x^c})}.
\end{align*}
These can and should be interpreted as a family of metrics on the local spin spaces at site $x \in \Z^d$ that measure how strongly a variation at site $y\in \Z^d$ can maximally change the difference in interaction energy between local spins $\varphi_x, \xi_x$. The criterion of the fineness of the coarse graining will involve the corresponding diameters of the preimages of the discrete variables. Recall that the diameter of a set $A$ with respect to some metric $d_{xy}$ is defined as usual as
\begin{align*}
    \text{diam}_{xy}(A) = \sup_{s,t \in A}d_{xy}(s,t). 
\end{align*}
Now we are ready to state the necessary criterion on the fineness of the discretisation. This is a very special case of the main result in \cite{van_enter_discrete_2011}. 

\begin{theorem}
    Suppose that the discretisation map $T:\mathbb{S}^1 \simeq [0,2\pi) \to \{1,\dots,q\}$ is such that the preimage of every $\varphi'\in\{1,\dots,q\}$ is an interval of the form $[\varphi'\lvert^l, \varphi'\lvert^r) \Subset [0,2\pi)$ 
    with left and right endpoints being given by 
    \begin{align*}
    &\varphi'\lvert^r = \sup\{ \varphi \in [0,2\pi): T(\varphi) = \varphi'\},
    \\\
    &\varphi'\lvert^l = \inf\{ \varphi \in [0,2\pi): T(\varphi) = \varphi'\}.
    \end{align*}
    and such that  
    \begin{align*}
        \sup_{x \in \Z^d}\sum_{y \neq x}\max_{\varphi'=1,\dots,q}\emph{diam}_{x,y}([\varphi'\lvert^l,\varphi'\lvert^r]) < 4,
    \end{align*}
    then we have $\overline{c}<1$,i.e., for every $\omega' \in \Omega'_q$ there exists a unique Gibbs measure $\mu_{\Z^d}[\omega']$ for the specification $\gamma^{\omega'}$.
\end{theorem}

Let us translate this into our situation and state a condition on the fineness of the discretisation, which we prove in Section~\ref{section:proofs}.
\begin{corollary}\label{corollary:choosing-q-sufficiently-large}
    For $q = q(\beta,\alpha) > \beta \pi \sum_{x \in \Z^d}\abs{x}^{-\alpha}$ we  have $\overline{c}<1$. 
\end{corollary}

So if $q$ is sufficiently large, for every fixed discrete configuration $\omega' \in \Omega'_q$ we can not only make sense of the specification $\gamma^{\omega'}$ but it actually gives rise to a \textit{unique} Gibbs measure on $\Omega$ which we will denote by $\mu_{\Z^d}[\omega'](d\omega)$. Note that one can actually construct this in a way such that this is a probability kernel from $\Omega'_q$ to $\Omega$, i.e., that it is also measurable as a function of the coarse-grained configuration, see \cite[Theorem~8.23]{georgii_gibbs_2011} or \cite[Section 2]{jahnel_class_2014}. 
From now on we will assume that $q$ is chosen sufficiently large and omit it from our notation, so we will for example just write $\Omega'$ instead of $\Omega'_q$. \medskip 

Before we can define the discrete specification $\gamma'$, there is one more object we need to introduce. 

\begin{lemma}\label{lemma:convergence-constrained-specifications}
    Let $\Lambda \Subset \Z^d$ and $((\gamma^{\omega'}\lvert_{\Lambda^c})_\Delta)_{\Delta\Subset \Lambda^c}$ the specification on $\Omega$ that one obtains as in Lemma~\ref{eqn:definition-constrained-specification} but by putting all potentials $\Phi_A$ with $A\cap \Lambda \neq \emptyset$ equal to $0$, i.e., we ignore all interactions with sites in $\Lambda$. Then, for every boundary condition $\omega \in \Omega$ with $\omega' = T(\omega)$, we have weak convergence 
    \begin{align*}
        (\gamma^{\omega'}\lvert_{ \Lambda^c})_\Delta(\cdot \lvert \omega) \rightharpoonup \mu_{\Lambda^c}[\omega'_{\Lambda^c}](\cdot),\qquad \text{as } \Delta\uparrow\Z^d,
    \end{align*}
    to the unique Gibbs measure $\mu_{\Lambda^c}[\omega'_{\Lambda^c}] \in \mathscr{G}(\gamma^{\omega'}\lvert_{\Lambda^c})$. Moreover, this convergence is uniform in $\omega$. 
\end{lemma}
In words, Lemma~\ref{lemma:convergence-constrained-specifications} establishes convergence of the specification conditioned on a coarse-grained configuration $\omega'_{\Lambda^c}$ inside any subvolume $\Lambda^c$ towards the corresponding limiting Gibbs measure. 

\subsubsection{The discrete specification and its regularity}

Now we are finally ready to write down the explicit definition of the discrete specification $\gamma'$. 
\begin{proposition}[Gibbsianness under discretisations]\label{lemma:explicit-form-discrete-specification}
    In the uniform Dobrushin regime, the discretisation image  $\mu'=\mu\circ T^{-1}$ of any continuous-spin Gibbs measure $\mu \in \mathscr{G}(\gamma)$ is a Gibbs measure for the specification $\gamma' = (\gamma'_\Delta)_{\Delta \Subset \Z^d}$, which is defined by  
    \begin{align*}
    \gamma'_\Lambda(\omega'_\Lambda \lvert \omega'_{\Lambda^c})
    =
    \frac{\mu_{\Lambda^c}[\omega'_{\Lambda^c}](\lambda^{\otimes \Lambda}(\exp(-\calH_\Lambda)\mathbf{1}_{[\omega'_\Lambda]}))}{\mu_{\Lambda^c}[\omega'_{\Lambda^c}](\lambda^{\otimes \Lambda}(\exp(-\mathcal{H}_\Lambda)))}.
\end{align*}
\end{proposition}

The proof is exactly the same as in \cite{jahnel_class_2014}, but we nevertheless include it below to make it more transparent why the $\gamma'_\Lambda$ are of the form above.  \medskip 

Via the Gibbs variational principle one can show that every translation-invariant discrete Gibbs measure has a continuous preimage and even give an explicit expression for the inverse in terms of a kernel. 

\begin{proposition}\label{proposition:kernel-correspondence-discrete-continuous}
    Let $\mu' \in \mathscr{G}(\gamma')$ be translation invariant. Then, its preimage w.r.t.~the discretisation map $T$ is given by 
    \begin{align*}
        \mu(d\omega) = \int \mu'(d\omega') \mu_{\Z^d}[\omega'](d\omega) \in \mathscr{G}(\gamma). 
    \end{align*}
    Moreover, the discretisation map $T$ induces a one-to-one map from $\emph{ex}\mathscr{G}_\theta(\gamma)$ to $\emph{ex}\mathscr{G}_\theta(\gamma')$. 
\end{proposition}

The proof of this makes use of the classical Gibbs variational principle see~\cite[Proposition 2.5]{jahnel_class_2014} for details.
Because every translation-invariant extremal Gibbs measure for $\gamma$ gets mapped to a translation-invariant extremal Gibbs measure for $\gamma'$, this in particular holds for the special family in Theorem \ref{theorem:continous-symmetry-breaking}. Let us parametrise this family by the angle $\varphi = \text{arg}(u)$ and write $\mu'_{\varphi}$ for the image of $\mu_\varphi$ under the discretisation map. 
\medskip 

As a last step, let us now show that the discretised speficiation depends smoothly on its boundary condition, in the sense of being of summable oscillations. We will need these regularity properties of $\gamma'$ to later apply the results from \cite{jahnel_dynamical_2023}. 

\begin{lemma}\label{lemma:regularity-discrete-specification}
    The discrete specification $\gamma'$ is translation invariant, nonnull, and satisfies 
    \begin{align*}
        \sum_{x \neq 0}\delta_x( \gamma'_0(\cdot)) < \infty.
    \end{align*}   
\end{lemma}

\subsection{Interacting particle systems}\label{section:interacting-particle-systems}

Now that we have introduced Gibbs measures, let us consider dynamics. We will study time-continuous, translation-invariant Markovian dynamics on $\Omega'$, namely interacting particle systems characterised by time-homogeneous generators $\mathscr{L}$ with domain $\text{dom}(\mathscr{L})$ and its associated Markovian semigroup $(S(t))_{t \geq 0}$. For interacting particle systems we adopt the notation and exposition of the standard reference \cite[Chapter 1]{liggett_interacting_2005}, but also refer the interested reader to the excellent exposition in \cite{swart_course_2022}. 

\subsubsection{Well-definedness via Liggett's criteria}
In our setting, the generator $\mathscr{L}$ is given via a collection of translation-invariant transition rates $c_\Delta(\eta, \xi_\Delta)$, in finite volumes $\Delta \Subset \Z^d$, which are continuous in the starting configuration $\eta \in \Omega$. 
These rates can be interpreted as the infinitesimal rate at which the particles inside $\Delta$ switch from the configuration $\eta_\Delta$ to $\xi_\Delta$, given that the rest of the system is currently in state $\eta_{\Delta^c}$. 
The full dynamics of the interacting particle system is then given as the superposition of these local dynamics, i.e., 
\begin{align*}
    \mathscr{L}f(\eta) = \sum_{\Delta \Subset \Z^d}\sum_{\xi_\Delta}c_\Delta(\eta, \xi_\Delta\eta_{\Delta^c})[f(\xi_\Delta \eta_{\Delta^c}) - f(\eta)].
\end{align*}
In \cite[Chapter 1]{liggett_interacting_2005} it is shown that the following two conditions are sufficient to guarantee the well-definedness for a given family of translation-invariant transition rates. 
\begin{enumerate}[\bfseries (L1)]
    \item The total rate at which the particle at a particular site changes its spin is uniformly bounded, i.e.,
    \begin{align*}
        \sum_{\Delta \ni 0} \sum_{\xi_{\Delta}}\norm{c_{\Delta}(\cdot, \xi_{\Delta}\cdot_{\Delta^c})}_{\infty} < \infty.
    \end{align*}
    \item The total influence of a single coordinate on all other coordinates is uniformly bounded, i.e.,
    \begin{align*}
        \sum_{\Delta \ni 0}\sum_{x \neq 0}\sum_{\xi_{\Delta}}\delta_x\left(c_{\Delta}(\cdot, \xi_{\Delta}\cdot_{\Delta^c})\right) < \infty.
    \end{align*} 
\end{enumerate}

Under these conditions, a core of the operator $\mathscr{L}$ is given by 
\begin{align*}
    D(\Omega) = \Big\{f \in C(\Omega): \ \vertiii{f} := \sum_{x \in \Z^d}\delta_x(f) < \infty \Big\}.
\end{align*}
Moreover, this core is invariant under the dynamics of the semigroup in the sense that for all $f \in D(\Omega)$ and $t\geq 0$ we have $S(t)f \in D(\Omega)$, see \cite[Theorem I.3.9]{liggett_interacting_2005}. 

Before we can introduce the family of transition rates we are interested in for our construction, we will need some more notation. For $\omega' \in \Omega'$ denote by $(\omega')^x$ the configuration which is identical to $\omega'$ everywhere except at site $x \in \Z^d$, where it is equal to $((\omega'(x) + 1)\text{ mod }q)$. Moreover, for a discrete spin variable $\varphi' \in \{1,\dots,q\}$ define 
\begin{align*}
    &\varphi'\lvert^r = \sup\{ \varphi \in [0,2\pi): T(\varphi) = \varphi'\},
    \\\
    &\varphi'\lvert^l = \inf\{ \varphi \in [0,2\pi): T(\varphi) = \varphi'\}.
\end{align*}
These are of course nothing but the right and the left endpoints of the arc of the circle (when parametrized by the angle) that gets mapped to the discrete variable $\varphi'$. 

\begin{lemma}\label{lemma:well-definedness-generator}
    The generator $\mathscr{L}$ with rates given by
    \begin{align}\label{eqn:transition-rates-irreversible}
        c(\omega', (\omega')^x)
        =
        \frac{\mu_{x^c}[\omega'_{x^c}](\exp(-\calH_x(\omega'_x\lvert^r, \cdot_{x^c})))}{\mu_{x^c}[\omega'_{x^c}](\lambda^x(\exp(-\calH_x) \mathbf{1}_{\omega'_x}))}
    \end{align}
    satisfies the well-definedness criteria $\mathbf{(L1)}$ and $\mathbf{(L2)}$, as well as the non-degeneracy condition $\mathbf{(ND)}$. 
\end{lemma}

Note that in contrast to the proof strategy in \cite{jahnel_class_2014} we will not need any further regularity of the rates like polynomial or even exponential decay of their oscillations.

\subsubsection{The infinitesimal rotation property for local functions}
Now, as a first step we show that – at least infinitesimally – the action of the semigroup $(S(t))_{t \geq 0}$ and the discretisation of the rotation 
\begin{align*}
R_t: \Omega_0 \to \Omega_0, \quad \omega \mapsto \left(\omega + t \mod 2 \pi\right), \quad t \geq 0, 
\end{align*}
agree on the extremal Gibbs measures when tested against local functions. Recall that we denote by $\mu'_\varphi$, $\varphi \in [0,2\pi)$, the Gibbs measure on $\Omega'$ which is the image of $\mu_u$ with $\text{arg}(u)=\varphi$ under the discretisation map, see Sections 3.2. and 3.3. for details, and note that the rotation map $R_t$ maps extremal Gibbs measures to extremal Gibbs measures in the continuous model. 

\begin{proposition}\label{proposition:local-infinitesimal-rotation-property}
    For all local functions $f$ and $t\geq 0$ it holds that
    \begin{align*}
        \frac{d}{d\varepsilon}\lvert_{\varepsilon=0}\mu'_{t+\varepsilon}(f) = \mu'_t(\mathscr{L}f). 
    \end{align*}
\end{proposition}

For the rates we are working with, one can show that the local functions are a core for the generator, see \cite[Theorem 4.30]{swart_course_2022}, but it is not closed under the dynamics of the semigroup. 

So we need to find a way to extend the infinitesimal rotation property from Proposition \ref{proposition:local-infinitesimal-rotation-property} to a larger class of functions which is closed under the dynamics of the semigroup.  

\subsubsection{Extension via local approximation}
The main technical helpers will be the following two approximation results.
The first one allows us to upgrade  convergence in the supremum-norm to convergence in the triple-norm, as long as we have some uniform control over the oscillations of the approximating functions. 
\begin{lemma}[Upgrading lemma]\label{lemma:upgrading-lemma}
    Let $f \in D(\Omega)$ and $(f_n)_{n \in \N}$ a sequence of functions with $f_n \in D(\Omega)$ for all $n \in \N$ such that $\norm{f-f_n}_\infty \to 0$ as $n \to \infty$ and
    \begin{align*}
        \sum_{x \in \Z^d}\sup_{n \in \N}\delta_x(f_n) < \infty. 
    \end{align*}
    Then we have $\vertiii{f-f_n}\to 0$ as $n \to \infty$. 
\end{lemma}

This we can now use to show that we can (usually) transfer infinitesimal properties of our dynamics from local functions to functions with summable oscillations. 
\begin{lemma}[Local approximation lemma]\label{lemma:local-approximation-lemma}
    Let $\mathscr{L}$ be the generator of the particle system of the form 
    \begin{align*}
        \mathscr{L}f(\eta) = \sum_{x \in \Z^d}c_x(\eta, \eta^x)[f(\eta^x)-f(\eta)].
    \end{align*}
    Moreover, assume that 
    \begin{align*}
        \mathbf{c} = \sup_{x \in \Z^d}\norm{c_x(\cdot, \cdot)}_\infty < \infty. 
    \end{align*}
    Let $f \in D(\Omega)$. Then, there exists a sequence of local functions $(f_n)_{n \in \N}$ such that $\vertiii{f_n - f} \to 0$ and $\norm{\mathscr{L}f - \mathscr{L}f_n}_\infty \to 0$ as $n \to \infty$. 
\end{lemma}

With these tools at hand, we can now show that the infinitesimal rotation property can be extended from the local functions to the whole space $D(\Omega)$. 

\begin{proposition}[Infinitesimal rotation property on $D(\Omega)$]\label{proposition:extending-infinitesimal-rotation-property}
    If $f \in D(\Omega)$, then for all $t\geq 0$
    \begin{align*}
        \frac{d}{d\varepsilon}\lvert_{\varepsilon=0}\mu'_{t+\varepsilon}(f) = \mu'_t(\mathscr{L}f). 
    \end{align*}
\end{proposition}

This is already very promising, because Proposition I.3.2 in \cite{liggett_interacting_2005} tells us that the space $D(\Omega)$ is closed under the action of the semigroup $(S(t))_{t \geq 0}$ generated by $\mathscr{L}$, or in other words, if $f \in D(\Omega)$, then $S(t)f \in D(\Omega)$ for all $t\geq 0$. Hence, we can apply the infinitesimal rotation property along trajectories of our semigroup. 

\subsubsection{A forward–backward argument to identify trajectories}
Fix a function $f\in D(\Omega)$ and $t > 0$. Consider the real-valued function 
\begin{align*}
    [0,t] \ni s\mapsto F(s) := \mu'_{s}(S(t-s)f). 
\end{align*}
So in some sense we are applying the deterministic rotation forward in time and our stochastic dynamics backward in time. 
Then (assuming that the derivative exists and is continuous) we can use the fundamental theorem of calculus to see that 
\begin{align*}
    \mu'_t(f) - \mu'_0(S(t)f) = F(t)-F(0) = \int_0^t \frac{d}{ds}\mu'_{s}(S(t-s)f) ds. 
\end{align*}
Hence, if we can show that 
\begin{align*}
    \frac{d}{ds}\mu'_{s}(S(t-s)f) = 0 
\end{align*}
for all $s\in (0,t)$, we have finally established that the rotation w.r.t.~the deterministic angle is reproduced by the stochastic evolution. 

\begin{proposition}\label{proposition:time-derivative-forward-backward}
    Let $f \in D(\Omega)$. Then, for all $s \in (0,t)$, the following limit exists and we have 
    \begin{align*}
        \lim_{\varepsilon\to 0}\frac{1}{\varepsilon}\left(\mu'_{s+\varepsilon}(S(t-s-\varepsilon)f)-\mu'_s(S(t-s)f)\right) =0. 
    \end{align*}
\end{proposition}

As mentioned above, the fundamental theorem of calculus now allows us to conclude the following. 

\begin{proposition}\label{proposition:full-rotation-property}
    For all $t,s \geq 0$ we have
    \begin{align*}
        \mu'_tS(s) = \mu'_{t+s \text{ mod } 2 \pi}.
    \end{align*}
\end{proposition}

This already establishes the existence of a well-defined interacting particle system that satisfies Part $(b)$ of Theorem \ref{theorem:main-result}. It remains to show the uniqueness of the equilibrium measure. 

\subsection{Uniqueness of the equilibrium measure}
By choosing the increasing sequence $T_n = 2\pi n$, $n \in \N$ in \cite[Proposition I.1.8(e)] {liggett_interacting_2005}, one directly sees that the probability measure defined by 
\begin{align*}
    \mu^* = \frac{1}{2\pi}\int_0^{2\pi}\mu'_t dt
\end{align*}
is time-stationary w.r.t.~the dynamics generated by $\mathscr{L}$. Moreover, since it is a convex combination of measures in $\mathscr{G}(\gamma')$, which is a convex set, it is also an element of $\mathscr{G}(\gamma')$. 
Hence, up to checking some regularity properties of the rates and the specification we are exactly in the situation of~\cite{jahnel_dynamical_2023}. More precisely, to show that $\mu^*$ is the only translation-invariant time-stationary measure for our dynamics, we first use that, by the results from \cite{jahnel_dynamical_2023}, every translation-invariant element of the attractor is necessarily also a Gibbs measure w.r.t.~the specification $\gamma'$. Recall that the attractor is defined by 
\begin{align*}
    \mathscr{A} = \Big\{ \nu \in \mathcal{M}_1(\Omega')\colon \exists \nu_0 \in \calM_1(\Omega) \text{ and } t_n \uparrow \infty \text{ s.t.} \lim_{n \to \infty}\nu_{t_n} = \nu\Big\},
\end{align*}
or in other words, the attractor $\mathscr{A}$ is the set of all weak limit points of the measure-valued dynamics induced by the Markov semigroup $(S(t))_{t\geq0}$.

\begin{proposition}\label{proposition:attractor-property}
    If $\nu'$ is translation invariant and in the attractor of the dynamics generated by $\mathscr{L}$, then it is necessarily a Gibbs measure w.r.t.~the specification $\gamma'$. 
    In particular, every translation-invariant and time-stationary measure is in $\mathscr{G}(\gamma')$. 
\end{proposition}
Because both our specification $\gamma'$ and the rates $c(\cdot, \cdot)$ are sufficiently regular, cf.~Lemma~\ref{lemma:regularity-discrete-specification} and Lemma \ref{lemma:well-definedness-generator}, this characterisation of the attractor follows directly from \cite[Theorem 10]{jahnel_dynamical_2023}. 
\medskip 

Now that we know that every time-stationary measure is also a Gibbs measure w.r.t.~the specification $\gamma'$, we can use the extremal decomposition and the rotation property to conclude that the equilibrium measure is indeed unique. 

\begin{proposition}[Uniqueness of the equilibrium measure]\label{proposition:uniqueness-equilibrium-measure}
    If $\nu'$ is a translation-invariant and time-stationary measure for the dynamics generated by $\mathscr{L}$, then 
    $\nu'=\mu^*$.
    In particular, there is a unique translation-invariant and time-stationary measure. 
\end{proposition}

This is shown exactly as in \cite{jahnel_class_2014}, but we nevertheless sketch the main idea in Section \ref{sec:proof-uniqueness-equilibrium-measure}. 
All in all, this last piece also establishes Part $(a)$ of Theorem \ref{theorem:main-result} and we conclude our construction. 

\begin{remark}\label{remark:stability-under-reversible-perturbations}
   The time-periodic orbit in Theorem \ref{theorem:main-result} is stable under reversible perturbations in the following sense.
   If $\mathscr{K}$ is the generator of an interacting particle system with translation-invariant rates that satisfy $\mathbf{(L1)}$ and $\mathbf{(L2)}$  and admits $\mu^*$ as a reversible measure, then the dynamics with generator $\mathscr{L}+\mathscr{K}$ is well-defined and still satisfies Properties $(a)$ and $(b)$ in Theorem \ref{theorem:main-result}, even if the semigroups generated by $\mathscr{L}$ and $\mathscr{K}$ do not commute. 
   To see this, first note that if $\mu^*$ is reversible for $\mathscr{K}$, then so are all other measures $\nu \in \mathscr{G}(\gamma')$. This allows us to get an analogue of Proposition \ref{proposition:local-infinitesimal-rotation-property} for the process generated by $\mathscr{K}+\mathscr{L}$. Now the whole procedure to extend the infinitesimal rotation property to $D(\Omega)$ and the forward-backward construction to conclude work just like for $\mathscr{L}$. 
\end{remark}

\newpage 

\section{Proofs}\label{section:proofs}

\subsection{Choosing the discretisation fine enough}

\begin{proof}[Proof of Corollary~\ref{corollary:choosing-q-sufficiently-large}]
    First note that by translation invariance it suffices to estimate the metrics $d_{0x}$ for $x \in \Z^d \setminus \{0\}$. 
    Recall that we performed an equidistant partition of the circle $\mathbb{S}^1 \simeq [0,2\pi)$ into $q$ arcs of length $2\pi/q$. By plugging in the definition of our Hamiltonian we get that, for two angles $\varphi_x, \xi_x \in [0,2\pi)$ that are in the same arc, we have 
    \begin{align*}
        d_{0,x}(\varphi_0, \xi_0)
        \leq 
        2\beta\abs{x}^{-\alpha}2\pi/q, 
    \end{align*}
    where we mainly used the Lipschitz-continuity of $u\mapsto \cos(u)$. Since we assume that $\alpha \in (d,2d)$, this is always summable, but in order to have
    \begin{align*}
        \sum_{x \neq 0}2\beta\abs{x}^{-\alpha}2\pi/q <4,  
    \end{align*}
    we need to choose 
    \begin{align*}
        q = q(\beta,\alpha) > \beta \pi \sum_{x \in \Z^d}\abs{x}^{-\alpha}.
    \end{align*}
    So one can also clearly see that this is an increasing function of the interaction strength, since it is increasing in $\beta$ and decreasing in $\alpha\in(d,2d)$. 
\end{proof}

\subsection{Convergence of the restricted and constrained specifications}

\begin{proof}[Proof of Lemma~\ref{lemma:convergence-constrained-specifications}]
    Since we have chosen our discretisation fine enough, we have $\overline{c}<1$, and, for fixed $\Lambda \Subset \Z^d$, the restricted and constrained specifications $\gamma^{\omega'}\lvert_{\Lambda^c}$ are in the Dobrushin uniqueness regime, uniform in the constraint $\omega'_\Lambda$ as well. Therefore, uniqueness holds by Theorem \ref{theorem:dobrushin-comparison} and the convergence is a direct consequence of \cite[Theorem 8.23]{georgii_gibbs_2011}. 
\end{proof}

\subsection{Explicit formula for the discrete specification}
\begin{proof}[Proof of Proposition~\ref{lemma:explicit-form-discrete-specification}]
    Let $\omega \in \Omega$ and $\omega'=T(\omega)$, then by definition of the discretisation
    \begin{align*}
        \mu'(\omega'_\Lambda) = \mu(\mathbf{1}_{[\omega'_\Lambda]}), \quad \Lambda \Subset \Z^d. 
    \end{align*}
    Further, by the martingale convergence theorem, the following convergence holds in $L^1(\mu)$ and almost surely as $\Lambda \uparrow \Z^d$
    \begin{align*}
        \mu'(\omega'_\Delta \lvert \omega'_{\Lambda \setminus \Delta})
        \to
        \mu(\mathbf{1}_{[\omega'_\Delta]}\lvert \mathcal{F}'_{\Delta^c})(\omega),
    \end{align*}
    where $\calF'_{\Delta^c}$ is the $\sigma$-algebra over $\Omega$ generated by the coarse-graining map $T$ only applied in $\Delta^c$. 

    Now note that, using the DLR equation for $\mu$, any of these conditional probabilities with finite-volume conditioning can be rewritten as 
    \begin{align}\label{eqn:intermediate-step}
        \mu'(\omega'_\Delta \lvert \omega'_{\Lambda \setminus \Delta})
        &=  \frac{\int\mu(d\omega)\gamma_\Lambda(\mathbf{1}_{[\omega'_\Delta]}\mathbf{1}_{[\omega'_{\Lambda \setminus \Delta}]}\lvert\omega_{\Lambda^c})}{\int \mu(d\omega)\gamma_\Lambda(\mathbf{1}_{[\omega'_{\Lambda \setminus \Delta}]}\lvert \omega_{\Lambda^c})}
        \\\
        &=
        \frac{\int\mu(d\omega)(\gamma^{\omega'}\lvert_{\Delta^c})_{\Lambda \setminus \Delta}(\lambda^{\otimes \Delta}(\exp(-\calH_\Delta)\mathbf{1}_{[\omega'_\Delta]})\lvert \omega_{\Lambda^c})}{\int\mu(d\omega)(\gamma^{\omega'}\lvert_{\Delta^c})_{\Lambda \setminus \Delta}(\lambda^{\otimes \Delta}(\exp(-\calH_\Delta))\lvert \omega_{\Lambda^c})} \nonumber 
    \end{align}
    where $\gamma^{\omega'}\lvert_{\Delta^c}$ is the specification on $\Omega$ that one obtains by putting all potentials $\Phi_A$ with $A\cap \Delta \neq \emptyset$ equal to $0$, i.e., we ignore all interactions with sites in $\Delta$, see Lemma~\ref{lemma:convergence-constrained-specifications}.
    Indeed, by writing out the definition of  $\gamma_\Lambda$ explicitly we see that 
    \begin{align*}
        Z_\Lambda(\omega_{\Lambda^c})
&\gamma_\Lambda(\mathbf{1}_{[\omega'_\Delta]}\mathbf{1}_{[\omega'_{\Lambda \setminus \Delta}]}\lvert \omega_{\Lambda^c}) 
        \\\
        &=
        \int \lambda^{\otimes \Lambda}(d\omega_{\Lambda})\mathbf{1}_{[\omega'_\Delta]}(\omega_{\Lambda})\mathbf{1}_{[\omega'_{\Lambda \setminus \Delta}]}(\omega_{\Lambda}) \exp(-\calH_\Lambda(\omega_{\Lambda}\omega_{\Lambda^c}))
        \\\
        &=
        \int \lambda^{\otimes \Lambda \setminus \Delta}(d\omega_{\Lambda \setminus \Delta})\mathbf{1}_{[\omega'_{\Lambda \setminus \Delta}]}(\omega_{\Lambda \setminus \Delta})\exp(-\calH_{\Lambda}^\Delta(\omega_{\Lambda \setminus \Delta}\omega_{\Lambda^c}))\times\\
        &\qquad\int \lambda^{\otimes \Delta}(d\omega_{\Delta})
        \mathbf{1}_{[\omega'_\Delta]}(\omega_{\Delta})\exp(-\calH_\Delta(\omega_{\Delta}\omega_{\Lambda\setminus\Delta}\omega_{\Lambda^c}))
        \\\
        &= 
        Z^{\omega'}_{\Lambda, \Delta}(\omega_{\Delta^c})
        \left(\gamma^{\omega'}\lvert_{\Delta^c}\right)_{\Lambda \setminus \Delta}(\lambda^{\otimes \Delta}(\exp(-\calH_\Delta)\mathbf{1}_{[\omega'_\Delta]}\lvert \omega_{\Lambda^c}),
    \end{align*}
    and similarly for the term in the denominator, where the prefactors 
    $Z$ (with appropriate indexes) are the  normalising factors as in \eqref{eqn:finite-volume-lebesgue-density}.
    Since these normalising factors appear both in the numerator and the denominator, they cancel out and one arrives at \eqref{eqn:intermediate-step}. 
    We smuggled the restricted and constrained specifications into the integrals, because we want to take the limit $\Lambda \uparrow \Z^d$ and the uniform Dobrushin uniqueness for these specifications gives us the pointwise convergence of these terms, and they are uniformly bounded because the interaction potential $\Phi$ is absolutely summable. More precisely, by Lemma~\ref{lemma:convergence-constrained-specifications}, we have for all quasilocal functions $f:\Omega \to \R$ that
    \begin{align*}
        (\gamma^{\omega'}\lvert_{\Delta^c})_{\Lambda \setminus \Delta}(f \lvert \omega_{\Lambda^c}) \to \mu_{\Delta^c}[\omega'_{\Delta^c}](f) 
    \end{align*}
    as $\Lambda \uparrow \Z^d$ uniformly in $\omega' = T(\omega)$. 
    Hence, we get via dominated convergence 
    \begin{align*}
        \lim_{\Lambda \uparrow \Z^d}&\frac{\int\mu(d\omega)(\gamma^{\omega'}\lvert_{\Delta^c})_{\Lambda \setminus \Delta}(\lambda^{\otimes \Delta}(\exp(-\calH_\Delta)\mathbf{1}_{[\omega'_\Delta]})\lvert \omega_{\Lambda^c})}{\int\mu(d\omega)(\gamma^{\omega'}\lvert_{\Delta^c})_{\Lambda \setminus \Delta}(\lambda^{\otimes \Delta}(\exp(-\calH_\Delta))\lvert \omega_{\Lambda^c})}
        \\\
        &=
        \frac{\mu_{\Delta^c}[\omega'_{\Delta^c}](\lambda^{\otimes \Delta}(\exp(-\calH_\Delta)\mathbf{1}_{[\omega'_\Delta]}))}{\mu_{\Delta^c}[\omega'_{\Delta^c}](\lambda^{\otimes \Delta}(\exp(-\calH_\Delta)))}
        =\gamma'_\Delta(\omega'_\Delta \lvert \omega'_{\Delta^c}).
    \end{align*}
    We hope that, by including this proof, it is a bit more clear why and how the Gibbs measures $\mu_{\Delta^c}[\omega'_{\Delta^c}]$ for the specification $\gamma^{\omega'}\lvert_{ \Delta^c}$ with open boundary conditions, i.e., no interactions with all sites in $\Delta$, show up in the definition of the discrete specification $\gamma'$. 
\end{proof}

\subsection{Regularity of the discrete specification}

To show that the discrete specification is smooth in the sense of having summable oscillations we will need the following regularity property of the interaction potential of the continuous spin system. 

\begin{lemma}\label{lemma:summable-oscillations-continuous-specification}
    Recall our interaction potential $\Phi=(\Phi_A)_{A \Subset \Z^d}$ given by 
    \begin{align*}
        \Phi_A(\omega)
        =
        \begin{cases}
            \abs{x-y}^{-\alpha} \omega_x \cdot \omega_y \quad &\text{if } A=\{x,y\}, \\\
            0 &\text{otherwise.}
        \end{cases}
    \end{align*}
    Then, we have 
    \begin{align*}
        \sup_{x \in \Z^d}\sum_{A \ni x}\vertiii{\Phi_A} < \infty. 
    \end{align*}
\end{lemma}

\begin{proof}
    First observe that, by translation invariance of the interaction potential, it suffices to show that
    \begin{align*}
        \sum_{x \in \Z^d}\vertiii{\Phi_{\{0,x\}}} < \infty.
    \end{align*}
    But then, for fixed $x \in \Z^d$, we have 
    \begin{align*}
        \forall y \in \Z^d \setminus \{0,x\}: \quad \delta_y \Phi_{\{0,x\}} = 0 
    \end{align*}
    and by symmetry also 
    \begin{align*}
        \delta_x \Phi_{\{0,x\}} = \delta_0 \Phi_{\{0,x\}}.
    \end{align*}
    For fixed $x \in \Z^d$, all $\omega \in \Omega$ and all $\varphi_1, \varphi_2 \in \mathbb{S}^1$ we have that
    \begin{align*}
        \abs{\Phi_{\{0,x\}}(\varphi_1 \omega_{x^c})-\Phi_{\{0,x\}}(\varphi_2\omega_{x^c})}
        \leq 
        \abs{x}^{-\alpha}. 
    \end{align*}
    By putting all of these ingredients together and recalling that $\alpha > d$ we see that 
    \begin{align*}
        \sum_{x \in \Z^d}\vertiii{\Phi_{\{0,x\}}}
        \leq 
        2 \sum_{x \in \Z^d}\abs{x}^{-\alpha} < \infty,
    \end{align*}
    as desired. 
\end{proof}

With this technical helper in place we are ready to show that the oscillations of the discrete specification are indeed summable. 

\begin{proof}[Proof of Lemma \ref{lemma:regularity-discrete-specification}]
    \textit{Translation invariance: } If one just considers the definition of the specification, translation invariance is not immediately clear. But, if $\mu$ is a translation-invariant Gibbs measure for the continuous-spin specification $\gamma$, then its image $\mu'$ under the discretisation map is also translation invariant. Now we know that $\gamma'$ is a regular version of the conditional probabilities of $\mu'$. But if $\mu'$ is translation invariant, then $\gamma'$ must also be translation invariant.  
    \medskip
    \newline 
    \textit{Nonnullness:} For all $\omega' \in \Omega'$ we have, by definition of $\gamma'$ and because the discretisation is finite and the continous potential $\Phi$ absolutely summable,  
    \begin{align*}
        \gamma_0'(\omega'_0 \lvert \omega'_{0^c})
        \geq 
        \exp(-2\norm{\calH_x}_\infty)2\pi/q > 0,
    \end{align*}
    so $\gamma'$ is indeed nonnull. 
    \medskip 
    \newline 
    \textit{Finite oscillations:} For $x \in \Z^d$ and $\omega', \eta' \in \Omega'$ with $\omega'\lvert_{x^c}= \eta'\lvert_{x^c}$ we have 
    \begin{align*}
        &\abs{\gamma'_0(\omega'_0 \lvert \omega'_{0^c}) - \gamma'_0(\eta'_0 \lvert \eta'_{0^c})}
        \\\
        &\quad\leq 
        \frac{q}{\pi}e^{2\norm{\calH_0}_\infty}\Big(\abs{\mu_{0^c}[\omega'_{0^c}](\lambda(\exp(-\calH_0)\mathbf{1}_{\omega'_0}))-\mu_{0^c}[\eta'_{0^c}](\lambda(\exp(-\calH_0)\mathbf{1}_{\omega'_0}))}
        \\\
        &\qquad+ 
        \abs{\mu_{0^c}[\omega'_{0^c}](\lambda(\exp(-\mathcal{H}_0)))-\mu_{0^c}[\eta'_{0^c}](\lambda(\exp(-\mathcal{H}_0)))}\Big),
    \end{align*}
    where we used the simple algebraic rule 
    \begin{align}\label{eq_alg}
        ad - bc = \frac{1}{2}[(a-b)(c+d)-(a+b)(c-d)]
    \end{align}
    and the basic estimate 
    \begin{align*}
        \mu_{0^c}[\omega'_{0^c}](\lambda(\exp(-\calH_0)\mathbf{1}_{\omega'_0}))\mu_{0^c}[\omega'_{0^c}](\lambda(\exp(-\mathcal{H}_0)))
        \geq 
        \exp(-2\norm{\calH_0}_\infty)2\pi/q > 0. 
    \end{align*}
    So in order to control the oscillations of the specification $\gamma'$, it suffices to control the difference of the integrals of certain quasilocal observables w.r.t.~the Gibbs measures $\mu_{0^c}[\omega'_{0^c}]$ and $\mu_{0^c}[\eta'_{0^c}]$. The main tool for this is Dobrushin's comparison estimate, Theorem~\ref{theorem:dobrushin-comparison}.
    For brevity, let us write  
    \begin{align*}
        \psi_1(\cdot) :=\lambda(\exp(-\calH_0(\cdot_0,\cdot_{0^c}))\mathbf{1}_{\omega'_0}), 
        \quad 
        \psi_2(\cdot) := \lambda(\exp(-\mathcal{H}_0(\cdot_0,\cdot_{0^c}))). 
    \end{align*}
    For any fixed discrete configuration $\omega' \in \Omega'$, the Gibbs measure $\mu_{0^c}[\omega'_{0^c}]$ is uniquely defined by the specification $\gamma^{\omega'}\lvert_{0^c}$. So if two discrete configurations $\omega'$ and $\eta'$ agree on $\Z^d \setminus \{x\}$, we have for all $y \in \Z^d \setminus 0$ and $\omega \in \Omega$
    \begin{align*}
        \norm{(\gamma^{\omega'}\lvert_{0^c})_y(\cdot \lvert \omega) - (\gamma^{\eta'}\lvert_{0^c})_y(\cdot \lvert \omega)}_{\text{TV}}\leq \mathbf{1}_x(y),
    \end{align*}
    because the kernels $\gamma_y^{\omega'}$ of the constrained specification only depend on the constraint $\omega'$ on the site $y \in \Z^d$, i.e., only on $\omega'_y$, see \eqref{eqn:definition-constrained-specification}. Combining this estimate with Dobrushin's comparison theorem yields for $\psi \in\{\psi_1, \psi_2\}$
    \begin{align*}
        \abs{\mu_{0 c}[\omega'_{0^c}](\psi) - \mu_{0^c}[\eta'_{0^c}](\psi)} 
        \leq 
        \sum_{y\neq0}\delta_y(\psi)D_{yx}(\gamma^{\omega'}\lvert_{0^c})
        \leq 
        \sum_{y \neq 0}\delta_y(\psi)\overline{D}_{yx},
    \end{align*}
    where we used that the constrained specifications $\gamma^{\omega'}$ are in the Dobrushin regime \textit{uniformly} in the constraint $\omega'$, since the discretisation is fine enough. In particular we have $\sum_{x \in \Z^d}\overline{D}_{yx}\leq C$ for some constant $C >0$ for all $y \in \Z^d$ and hence 
    \begin{align*}
        \sum_{x \neq 0}\sum_{y\neq 0}\delta_y(\psi)\overline{D}_{yx} 
        \leq 
        C\sum_{y \in \Z^d}\delta_y(\psi). 
    \end{align*}
    Now, for every fixed $y \in \Z^d$, we can use the elementary inequality
    \begin{align}\label{eq_el}
        \abs{e^s - e^t} \leq \abs{s-t}e^{\max(\abs{s},\abs{t})}, \quad s,t \in \R, 
    \end{align}
    to obtain
    \begin{align*}
        \sum_{y \in \Z^d}\delta_y(\psi)
        \leq 
        e^{\norm{\calH_0}_\infty}\sum_{y \in \Z^d}\sum_{A \ni 0}\delta_y(\Phi_A) 
        =
        e^{\norm{\calH_0}_\infty}\sum_{A \ni 0}\vertiii{\Phi_A} < \infty.
    \end{align*}
    To see that the right-hand side in the last line is actually finite we used Lemma~\ref{lemma:summable-oscillations-continuous-specification}. 
    Since the right-hand side is uniform for all $\omega',\eta'$ with $\omega'\lvert_{x^c}= \eta'\lvert_{x^c}$, we get that 
    \begin{align*}
        \sum_{x \neq 0}\delta_x(\gamma'_0(\cdot)) < \infty,
    \end{align*}
    which, by translation invariance of $\gamma'$, is exactly what we want. 
\end{proof}

\subsection{Well-definedness of the generator}

For the proof we will again use the estimate on the oscillations of the continuous interaction potential $\Phi$ from Lemma~\ref{lemma:summable-oscillations-continuous-specification}.  
\begin{proof}[Proof of Lemma \ref{lemma:well-definedness-generator}]
    \textit{Ad $\mathbf{(L1)}$:} We only need to check if for the rate at which a single particle rotates the spin is bounded. By definition of the rates we have, for any discrete configuration $\omega' \in \Omega'$ and $x \in \Z^d$,
    \begin{align*}
        \abs{c(\omega',(\omega')^x)} \leq \exp(2\norm{\calH_x}_\infty)q/(2\pi) < \infty. 
    \end{align*}
    So the rate at which at particle at a particular site changes its spin is uniformly bounded. 
    \medskip
    \newline 
    \textit{Ad $\mathbf{(L2)}$:} 
    By using again the simple algebraic rule~\eqref{eq_alg} 
    and the basic estimate 
    \begin{align}\label{eqn:basic-estimate}
        \mu_{0^c}[\omega'_{0^c}](\exp(-\calH_0(\omega'_0\lvert^r, \cdot_{0^c}))) \cdot \mu_{0^c}[\omega'_{0^c}](\lambda(\exp(-\calH_0) \mathbf{1}_{\omega'_0}))
        \geq 
        \tfrac{2\pi}{q}\exp(-2\norm{\calH_0}),
    \end{align}
    we obtain for $\omega', \eta' \in \Omega'$ with $\omega'\lvert_{x^c} = \eta'\lvert_{x^c}$
    \begin{align*}
        &\abs{c(\omega',(\omega')^0)-c(\eta',(\eta')^0)}
        \\\
        &\leq 
        \tfrac{q}{\pi}e^{2\norm{\calH}_\infty}\Big(
        \abs{\mu_{0^c}[\omega'_{0^c}](\exp(-\calH_0(\omega'_0\lvert^r, \cdot_{0^c})))-\mu_{0^c}[\eta'_{0^c}](\exp(-\calH_0(\omega'_0\lvert^r, \cdot_{0^c})))}
        \\\
        &\qquad+
        \abs{\mu_{0^c}[\omega'_{0^c}](\lambda(\exp(-\calH_0(\cdot, \cdot_{0^c})\mathbf{1}_{\omega'_0}))-\mu_{0^c}[\eta'_{0^c}](\lambda(\exp(-\calH_0(\cdot, \cdot_{0^c})\mathbf{1}_{\omega'_0}))}      \Big), 
    \end{align*}
    and therefore it suffices to control the difference of the integrals of certain quasilocal observables with respect to the Gibbs measures $\mu_{0^c}[\omega'_{0^c}]$ and $\mu_{0^c}[\eta'_{0^c}]$. Hence, we can use precisely the same arguments as in the proof for finite oscillations in Lemma~\ref{lemma:regularity-discrete-specification} above to see that      
    \begin{align*}
        \sum_{x \neq 0}\delta_x(c(\cdot, \cdot^0)) < \infty,
    \end{align*}
   which, for translation-invariant rates, is exactly the condition $\mathbf{(L2)}$. 
    \medskip 

    \textit{Ad $\mathbf{(ND)}$:} Let $\Delta \Subset \Z^d$ and $\eta, \xi$ agree on $\Delta^c$. Then we can simply proceed by using some enumeration of the sites of $\Delta = \{x_1, \dots, x_n\}$ and then – one after the other – rotate the spin at each of the sites one discrete step at a time. The transition rate for every step is strictly positive as one can see by considering \eqref{eqn:basic-estimate}. 
\end{proof}

\subsection{The infinitesimal-rotation property for local functions}
For the proof of Proposition~\ref{proposition:local-infinitesimal-rotation-property} we will make use of the following identity for the constrained and restricted Gibbs measures. 

\begin{lemma}\label{lemma:constrained-gibbs-measures-useful-identity}
    For all $\Lambda \Subset \Z^d$ finite, and $\calF_{\Lambda^c}$-measurable $\varphi:\Omega \to \R$ it holds that 
    \begin{align*}
        \mu_{\Z^d}[\omega'](\varphi(\cdot_{\Lambda^c}))
        =
        \frac{\mu_{\Lambda^c}[\omega'_{\Lambda^c}](\varphi(\cdot_{\Lambda^c})\lambda^{\otimes \Lambda}(\exp(-\calH_\Lambda)\mathbf{1}_{[\omega'_\Lambda]}))}{\mu_{\Lambda^c}[\omega'_{\Lambda^c}](\lambda^{\otimes \Lambda}(\exp(-\calH_\Lambda) \mathbf{1}_{[\omega'_\Lambda}])}
    \end{align*}
\end{lemma}

\begin{proof}
    We use the construction of these measures as limits of the conditioned continuous-spin specifications $\gamma^{\omega'}$. 
    In the uniform Dobrushin uniqueness regime, see  \cite[Theorem 8.23]{georgii_gibbs_2011}, we can use that, for every $\omega \in \Omega$, the left-hand side can be written as
    \begin{align*}
        \mu_{\Z^d}[\omega'](\varphi(\cdot_{\Lambda^c}))
        =
        \lim_{\Delta \uparrow \Z^d}\gamma_\Delta^{\omega'}(\varphi(\cdot_{\Lambda^c})\lvert \omega). 
    \end{align*}
    For fixed $\Delta \Subset \Z^d$ such that $\Lambda \Subset \Delta$, we can write out the terms explicitly as
    \begin{align*}
        &\gamma_\Delta^{\omega'}(\varphi(\cdot_{\Lambda^c})\lvert \omega) 
        \\\
        &=
        \tfrac{1}{Z_{\Delta}^{\omega'}(\omega_{\Delta^c})}
        \int_{[\omega'_\Delta]} \lambda^{\otimes \Delta}(d\omega_\Delta)e^{-\calH_\Delta(\omega_{\Delta}\omega_{\Delta^c})}\varphi(\omega_{\Lambda^c})
        \\\
        &=
        \tfrac{1}{Z_{\Delta}^{\omega'}(\omega_{\Delta^c})}\int_{[\omega'_{\Delta \setminus \Lambda}]}\lambda^{\otimes \Delta \setminus \Lambda}(d\omega_{\Delta \setminus \Lambda})\int_{[\omega'_\Lambda]} \lambda^{\otimes \Lambda}(d\omega_\Lambda)
        e^{-\calH_\Delta^\Lambda(\omega_{\Delta} \omega_{\Delta^c})}e^{-\calH_\Lambda(\omega_\Delta \omega_{\Delta^c})}\varphi(\omega_{\Lambda^c}), 
    \end{align*}
    where 
    \begin{align*}
        \calH_\Delta^\Lambda(\omega) = \sum_{A\colon A \cap \Delta \neq \emptyset, A \cap \Lambda = \emptyset}\Phi_A(\omega). 
    \end{align*}
    The above can be rewritten as 
    \begin{align*}
        &\frac{1}{Z_{\Delta}^{\omega'}(\omega_{\Delta^c})}
        \int_{[\omega'_{\Delta \setminus \Lambda}]}\lambda^{\otimes \Delta \setminus \Lambda}(d\omega_{\Delta \setminus \Lambda})
        e^{-\calH_\Delta^\Lambda(\omega_{\Lambda^c})}
        \varphi(\omega_{\Lambda^c})
        \int_{[\omega'_\Lambda]} \lambda^{\otimes \Lambda}(d\omega_\Lambda)
        e^{-\calH_\Lambda(\omega_\Delta \omega_{\Delta^c})}.
    \end{align*}
    Further,  note that the normalisation constants of $\gamma_\Delta^{\omega'}$ and $(\gamma^{\omega'}\lvert_{\Lambda^c})_\Delta$ satisfy
    \begin{align*}
        Z_{\Delta}^{\omega'}(\omega_{\Delta^c}) 
        =
        Z_{\Delta, \Lambda}^{\omega'}(\omega_{\Delta^c})(\gamma^{\omega'}\lvert_{\Lambda^c})_\Delta(\exp(-\calH_\Lambda)\mathbf{1}_{[\omega_\Lambda]}\lvert \omega_{\Delta^c}). 
    \end{align*}
    Indeed, by expanding the left-hand side and refactoring some of the terms one gets  
    \begin{align*}
        Z_\Delta^{\omega'}(\omega_{\Delta^c}) 
        &=
        \int_{[\omega'_\Delta]} \lambda^{\otimes\Delta}(d\omega_\Delta) e^{-\calH_\Delta(\omega_\Delta \omega_{\Delta^c})}
        \\\
        &=
        \int_{[\omega'_{\Delta \setminus \Lambda}]} \lambda^{\otimes \Delta \setminus \Lambda}(d\omega_{\Delta \setminus \Lambda})e^{-\calH_\Delta^\Lambda(\omega_{\Delta\setminus\Lambda} \omega_{\Delta^c})}
        \int_{[\omega'_\Lambda]}\lambda^{\otimes\Lambda}(d\omega_\Lambda)e^{-\calH_\Lambda(\omega_\Lambda\omega_{\Delta\setminus\Lambda}\omega_{\Delta^c})}
        \\\
        &= 
        Z_{\Delta, \Lambda}^{\omega'}(\omega_{\Delta^c})(\gamma^{\omega'}\lvert_{\Lambda^c})_\Delta(e^{-\calH_\Lambda}\mathbf{1}_{[\omega'_\Lambda]}\lvert \omega_{\Delta^c}). 
    \end{align*}
    By plugging this in we obtain
    \begin{align*}
\gamma^{\omega'}_\Delta(\varphi(\cdot_{\Lambda^c})\lvert \omega) 
        =
        \frac{(\gamma^{\omega'}\lvert_{\Lambda^c})_\Delta(\varphi(\cdot_{\Lambda^c})\lambda^{\otimes \Lambda}(\exp(-\calH_\Lambda)\mathbf{1}_{[\omega'_\Lambda]})\lvert \omega_{\Delta^c})}{(\gamma^{\omega'}\lvert_{\Lambda^c})_\Delta(\exp(-\calH_\Lambda)\mathbf{1}_{[\omega'_\Lambda]}\lvert \omega_{\Delta^c})}.
    \end{align*}
    Now we can take the limit $\Delta \uparrow \Z^d$ and use that the specifications $\gamma^{\omega'}$ are in the Dobrushin uniqueness regime to get the claimed identity. 
\end{proof}

Before we can start with the proof of Proposition~\ref{proposition:local-infinitesimal-rotation-property} we need one more Lipschitz-type estimate for the action of the deterministic rotation on the discretised Gibbs measures. 

\begin{lemma}\label{lemma:probability-of-jumping}
    Consider the continuum Gibbs measures $\{\mu_t: t \in [0,2\pi)\} \subset \emph{ex}\mathscr{G}_\theta(\Phi)$. Then, there exists a constant $C>0$ such that 
    \begin{align*}
        \forall \varepsilon > 0 \ \forall j \in \Z^d \ \forall t \geq 0: 
        \quad 
        \mu_t[\omega: T(\omega_j - \varepsilon) = T(\omega_j) -1] \leq C \varepsilon. 
    \end{align*}
\end{lemma}

\begin{proof}
    We will use the short-hand notation $A_j = \{\omega: T(\omega_j - \varepsilon) = T(\omega) -1\}$. By the DLR equation we can estimate 
    \begin{align*}
        \mu_t[A_j]
        \leq
        \sup_{\eta \in \Omega}\gamma_j(A_j \lvert \eta) 
        \leq 
        \frac{\varepsilon q e^{\norm{\calH_j}_\infty}}{2\pi e^{-\norm{\calH_j}_\infty}} = C\varepsilon. 
    \end{align*}
     Note that by the translation invariance of the Hamiltonian $\calH$, the right-hand side does not depend on the site $j \in \Z^d$ (and since we take the supremum also not on the angle $t$). 
\end{proof}

This can be used to get the following quantitative estimate on taking expectations w.r.t.~slightly rotated discretised Gibbs measures $\mu'_{t+\varepsilon}$ and $\mu'_t$. 

\begin{lemma}\label{lemma:lipschitz-continuity-of-rotation}
    There exists a uniform constant $C>0$ such that for all $g \in D(\Omega)$ and $t\geq 0$ we have 
    \begin{align*}
        \forall \varepsilon>0: \quad \abs{\mu'_{t+\varepsilon}(g)-\mu'_t(g)} \leq C \vertiii{g}\varepsilon. 
    \end{align*}
\end{lemma}

\begin{proof}
    By definition of the discretised measures and telescoping we can write 
    \begin{align*}
        \mu'_{t+\varepsilon}(g) - \mu'_t(g)
        &=
        \int \mu_t(d\omega)\left(g(T(\omega-\varepsilon\mathbf{1}_{\Z^d})-g(\omega)\right)
        \\\
        &=
        \int \mu_t(d\omega) \sum_{j \in \Z^d}\left(g(T(\omega-\varepsilon\mathbf{1}_{[0,\dots,\sigma(j)]})-g(T(\omega - \mathbf{1}_{[0,\dots, \sigma(j)-1]})\right)
        \\\
        &\leq 
        \sum_{j \in \Z^d}\delta_j(g)\mu_t[\omega: \ T(\omega_j - \varepsilon) = T(\omega) -1],
    \end{align*}
    where we used some arbitrary ordering $\sigma:\Z^d \to \N$ for the telescoping. This we can now estimate via Lemma \ref{lemma:probability-of-jumping} to obtain 
    \begin{align*}
        \abs{\mu'_{t+\varepsilon}(g) - \mu'_t(g)} \leq C\varepsilon\sum_{j \in \Z^d}\delta_j(g) = C\varepsilon\vertiii{g},
    \end{align*}
    as desired.
\end{proof}

After establishing these technical estimates, we can now show that the discretisation of the deterministic rotation and the action of the semigroup agree infinitesimally, at least when tested against local functions. 

\begin{proof}[Proof of Proposition \ref{proposition:local-infinitesimal-rotation-property}]
    Since the coarse graining is finite, i.e., $\abs{\Omega'_0} = q < \infty$, we can assume without loss of generality that $f=\mathbf{1}_\eta$ for some finite volume $\Lambda \Subset \Z^d$ and a finite-volume configuration $\eta' \in \Omega'_\Lambda$.
    Write $\rho_\Lambda = d\gamma_\Lambda/d\lambda^{\otimes \Lambda}$ for the Lebesgue density of the local specification in $\Lambda$, see \eqref{eqn:finite-volume-lebesgue-density}.
    By our Lipschitz estimate in Lemma~\ref{lemma:lipschitz-continuity-of-rotation}, the following derivative exists almost everywhere and we can use the DLR equation to calculate 
    \begin{align*}
        \frac{d}{d\varepsilon}\lvert_{\varepsilon=0}\mu'_{t+\varepsilon}(\mathbf{1}_{\eta'})
    &=
        \int \mu_t(d\omega)\frac{d}{d\varepsilon}\lvert_{\varepsilon=0} \left(\prod_{i \in \Lambda}\int_{\eta'_i\lvert^l - \varepsilon}^{\eta'_i\lvert^r- \varepsilon}\right) \lambda^{\otimes \Lambda}(d\varphi_\Lambda)\rho_\Lambda(\varphi_\Lambda\omega_{\Lambda^c})
        \\\
        &=
        \sum_{j \in \Lambda}\int \mu_t(d\omega) \left(\prod_{i \in \Lambda \setminus \{j\}}\int_{\eta'_i\lvert^l}^{\eta'_i\lvert^r}\right)\lambda^{\otimes \Lambda \setminus\{j\}}(d\varphi_{\Lambda \setminus \{j\}})\times\\
        &\qquad \qquad \qquad
        \left(\rho_\Lambda(\eta'_j\lvert^l \varphi_{\Lambda \setminus \{j\}} \omega_{\Lambda^c})-\rho_\Lambda(\eta'_j\lvert^r\varphi_{\Lambda \setminus \{j\}} \omega_{\Lambda^c})\right). 
    \end{align*}
    On the other hand, we know that, by definition of the generator,
    \begin{align}\label{eqn:generator-calculation}
        \mu'_t(\mathscr{L}\mathbf{1}_\eta) 
        &= 
        \sum_{j \in \Lambda}\Big(\int_{\{(\omega')^j_\Lambda = \eta\}}\hspace{-0.1cm}\mu'_t(d\omega')c(\omega',(\omega')^j)-\int_{\{\omega'_\Lambda = \eta\}}\mu'_t(d\omega')c(\omega', (\omega')^j)\Big).
    \end{align}
    It remains to show that we can rewrite these terms appropriately.
    For fixed $j\in \Lambda$ we can plug in the definition of the transition rates to get
    \begin{align*}
        &\int\mu'_t(d\omega')\mathbf{1}_{\{\omega'_\Lambda = \eta\}}(\omega')c(\omega', (\omega')^j)
        \\\
        &=\int \mu'_t(d\omega')\mathbf{1}_{\{\omega'_\Lambda = \eta\}}(\omega')\frac{\mu_{\Lambda^c}[\omega'_{\Lambda^c}](\lambda^{\otimes \Lambda\setminus \{j\}}(\exp(-\calH_\Lambda(\eta'_j\lvert^r, \cdot_{\Lambda \setminus \{j\}}, \cdot_{\Lambda^c}))))}{\mu_{\Lambda^c}[\omega'_{\Lambda^c}](\lambda^{\otimes \Lambda}(\exp(-\calH_\Lambda(\cdot_\Lambda,\cdot_\Lambda^c))))}
        \\\
        &=\int \mu'_t(d\omega')\mathbf{1}_{\{\omega'_\Lambda = \eta\}}(\omega')\mu_{\Z^d}[\omega']\left(\frac{\lambda^{\otimes \Lambda\setminus \{j\}}(\exp(-\calH_\Lambda(\eta'_j\lvert^r, \cdot_{\Lambda \setminus \{j\}}, \cdot_{\Lambda^c}))))}{\lambda^{\otimes \Lambda}(\exp(-\calH_\Lambda(\cdot_\Lambda,\cdot_{\Lambda^c})))}\right), 
    \end{align*}
    where, in the last step, we used Lemma \ref{lemma:constrained-gibbs-measures-useful-identity} with \begin{align*}
        \varphi(\omega) 
        =
        \frac{\lambda^{\otimes \Lambda\setminus \{j\}}(\exp(-\calH_\Lambda(\omega'_j\lvert^r, \cdot_{\Lambda \setminus \{j\}}, \cdot_{\Lambda^c}))))}{\lambda^{\otimes \Lambda}(\exp(-\calH_\Lambda(\cdot_\Lambda,\cdot_{\Lambda^c})))}, \quad \omega \in \Omega. 
    \end{align*}
Now, let us rewrite the integrand to see that 
    \begin{align*}
        &\int \mu'_t(d\omega')\mathbf{1}_{\{\omega'_\Lambda = \eta\}}(\omega')\mu_{\Z^d}[\omega']\left(\frac{\lambda^{\otimes \Lambda\setminus \{j\}}(\exp(-\calH_\Lambda(\eta'_j\lvert^r, \cdot_{\Lambda \setminus \{j\}}, \cdot_{\Lambda^c}))))}{\lambda^{\otimes \Lambda}(\exp(-\calH_\Lambda(\cdot_\Lambda,\cdot_{\Lambda^c})))}\right)
        \\\
        &=\int \mu'_t(d\omega')\mu_{\Z^d}[\omega']
        \left( 
            \frac{\mathbf{1}_{[\eta']}(\cdot)}{\gamma_\Lambda(\mathbf{1}_{[\eta']}\lvert \cdot)}
            \left(\prod_{i \in \Lambda \setminus \{j\}}\int_{\eta_i\lvert^l}^{\eta_i\lvert^r} \right)
            \lambda^{\otimes \Lambda \setminus \{j\}}(d\varphi_{\Lambda \setminus \{j\}})\rho_\Lambda(\eta'_j\lvert^r \varphi_{\Lambda\setminus \{j\}} \cdot_{\Lambda^c})
        \right) 
        \\\
        &=\int \mu_t(d\omega)\left(\prod_{i \in \Lambda \setminus \{j\}}\int_{\eta'_i\lvert^l}^{\eta'_i\lvert^r}\right)\lambda^{\otimes \Lambda \setminus \{j\}}(d\varphi_{\Lambda \setminus \{j\}})\rho_\Lambda(\eta'_j\lvert^r \varphi_{\Lambda \setminus \{j\}} \omega_{\Lambda^c}), 
    \end{align*}
    where we used the DLR equation and the fact that the kernel $\mu_{\Z^d}[\omega'](d\omega)$ gives us the correspondence between continuous and discrete Gibbs measures, see Proposition~\ref{proposition:kernel-correspondence-discrete-continuous}. 
    By proceeding similarly with the other terms in~\eqref{eqn:generator-calculation} and putting the two back together again we obtain the claimed identity for all local functions $f:\Omega' \to \R$. 
\end{proof}

\subsection{Extending the infinitesimal rotation property}
We first prove our technical helpers.

\begin{proof}[Proof of Lemma \ref{lemma:upgrading-lemma} -- the upgrading lemma]
    Let $A \Subset \Z^d$ be an arbitrary finite subset of $\Z^d$. Then, we can write 
    \begin{align*}
        \vertiii{f-f_n}
        &=
        \sum_{x \in \Z^d}\sup_{\eta \in \Omega,\ k=1,\dots,q-1}\abs{f(\eta^{x,k})-f_n(\eta^{x,k})-f(\eta)+f_n(\eta)},
    \end{align*}
    where $\eta^{x,k}$ is the configuration that is equal to $\eta$ everywhere except at the site $x \in \Z^d$, where the spin is rotated by $k$ discrete steps in clockwise direction, modulo $q$.
    Now, there are two different ways in which we can estimate the terms in this sum. On the one hand,
    \begin{align*}
        \abs{f_n(\eta^{x,k})-f(\eta^{x,k})-f_n(\eta)+f(\eta)}
        \leq 
        \abs{f_n(\eta^{x,k})-f_n(\eta)} + \abs{f(\eta^{x,k})-f(\eta)}
        \leq 
        2 \delta_x(f)
    \end{align*}
    and on the other hand 
    \begin{align*}
        \abs{f_n(\eta^{x,k})-f(\eta^{x,k})-f_n(\eta)+f(\eta)}
        &\leq 
        \abs{f_n(\eta^{x,k})-f(\eta^{x,k})} + \abs{f_n(\eta)-f(\eta)}
        \\\
        &\leq 
        2\norm{f_n - f}_\infty. 
    \end{align*}
    By using the former when $x \notin A$ and the latter when $x \in A$ we get that 
    \begin{align*}
        \vertiii{f_n - f}
        &\leq 
        2\sum_{x\in A}\norm{f - f_n}_\infty + \sum_{x \notin A}\left(\delta_x(f) + \delta_x(f_n)\right)
        \\\
        &=
        2\abs{A}\norm{f-f_n}_\infty + \sum_{x \notin A}\left(\delta_x(f) + \delta_x(f_n)\right). 
    \end{align*}
    By letting $n$ tend to infinity and using that $\norm{f_n-f}_\infty \to 0$ this inequality implies that
    \begin{align*}
        \limsup_{n \to \infty}\vertiii{f-f_n} \leq \sum_{x \notin A}\left(\delta_x(f)+\sup_{n \in \N}\delta_x(f_n)\right).
    \end{align*}
    Now we can let $A \uparrow \Z^d$ to see that $\vertiii{f-f_n}$ goes to zero as $n$ tends to infinity. 
\end{proof}

Now we can use the upgrading lemma to show our local approximation result. 
\begin{proof}[Proof of \ref{lemma:local-approximation-lemma} -- the local approximation lemma]
    Choose a sequence $\Lambda_n \uparrow \Z^d$ of finite sets, fix some reference configuration $\omega \in \Omega$ and put 
    \begin{align*}
        f_n(\eta) := f(\eta_{\Lambda_n}\omega_{\Lambda_n^c}). 
    \end{align*}
    Then, each $f_n$ is clearly a local function, since it only depends on the coordinates in $\Lambda_n$. In particular, we have $f_n \in D(\Omega)$ for all $n \in \N$. By uniform continuity of $f$ and the definition of the $f_n$ we directly get that $\norm{f-f_n}_\infty \to 0$ as $n$ tends to infinity and for every site $x \in \Z^d$ we have 
    \begin{align*}
        \sup_{n \in \N}\delta_x(f_n) 
        \leq 
        \delta_x(f). 
    \end{align*}
    Hence, by the upgrading lemma this implies $\vertiii{f_n-f} \to 0$ as $n$ tends to infinity. 
    Moreover, we have for any $\eta \in \Omega$
    \begin{align*}
        \abs{\mathscr{L}f_n(\eta) - \mathscr{L}f(\eta)}
        \leq 
        \sum_{x \in \Z^d}\norm{c_x(\cdot, \cdot)}_\infty \abs{f_n(\eta^x)-f_n(\eta) - f(\eta^x)+f(\eta)}
        \leq 
        2 \mathbf{c} \vertiii{f_n -f}, 
    \end{align*}
     and we have already shown that the right-hand side converges to zero as $n$ tends to infinity. 
\end{proof}

We now extend the infinitesimal rotation property to $D(\Omega)$. For this step, we will use the following corollary to the upgrading lemma. 

\begin{lemma}
    Let $f \in D(\Omega)$, and $(S(t))_{t \geq 0}$ the semigroup of an interacting particle system with generator $\mathscr{L}$ such that its rates satisfy Conditions $\mathbf{(L1)}$ and $\mathbf{(L2)}$. Then, we have that
    \begin{align*}
        \lim_{t \downarrow 0}\vertiii{S(t)f-f} = 0.
    \end{align*}
\end{lemma}

\begin{proof}
    Fix $f \in D(\Omega)$. Then, from \cite[Theorem I.3.9(c)]{liggett_interacting_2005}, we know that there exists a constant $C=C(T)>0$ such that 
    \begin{align*}
        \forall x \in \Z^d: 
        \quad 
        \sup_{t \in [0,T]}\delta_x(S(t)f) \leq C \delta_x(f), 
    \end{align*}
    and by strong-continuity of the semigroup we know that $\norm{S(t)f-f}_\infty \to 0$ as $t \to 0$. So we can apply the upgrading lemma and get $\vertiii{S(t)f-f} \to 0$ as $t \to 0$. 
\end{proof}

With these approximation tools at hand, we can now finally show that the infinitesimal-rotation property actually holds for all functions $f\in D(\Omega)$ and not just for local functions. 
\begin{proof}[Proof of Proposition \ref{proposition:extending-infinitesimal-rotation-property}]
    Let $f\in D(\Omega)$. Then, by the local approximation Lemma~\ref{lemma:local-approximation-lemma}, there exists a sequence of local functions $(f_n)_{n \in \N}$ such that both $\vertiii{f_n -f} \to 0$ and $\norm{\mathscr{L}f - \mathscr{L}f_n}_\infty \to 0$ as $n \to \infty$. 
    We can apply the Lipschitz estimate from Lemma \ref{lemma:lipschitz-continuity-of-rotation} with $g=f-f_n$ to see that there exists a constant $C$ that does not depend on $n$ or $\varepsilon$ such that  
    \begin{align*}
        \tfrac{1}{\varepsilon}\abs{\mu'_{t+\varepsilon}(f-f_n)-\mu'_t(f -f_n)}
        \leq 
        C \vertiii{f-f_n}.
    \end{align*}
    By combining these two estimates we get 
    \begin{align*}
        &\abs{\tfrac{1}{\varepsilon}\left(\mu'_{t+\varepsilon}(f)-\mu'_t(f)\right) - \mu'_t(\mathscr{L}f)}
        \\\
        &\qquad\leq 
        \tfrac{1}{\varepsilon}\abs{\mu'_{t+\varepsilon}(f-f_n)-\mu'_t(f-f_n)}
        +
        \abs{\mu'_t(\mathscr{L}f-\mathscr{L}f_n)}\\
        &\qquad\qquad
        +
        \abs{\tfrac{1}{\varepsilon}\left(\mu'_{t+\varepsilon}(f_n)-\mu'_t(f_n)\right) - \mu'_t(\mathscr{L}f_n)}
        \\\
        &\qquad\leq 
         C \vertiii{f_n -f} + \norm{\mathscr{L}f-\mathscr{L}f_n}_\infty +
        \abs{\tfrac{1}{\varepsilon}\left(\mu'_{t+\varepsilon}(f_n)-\mu'_t(f_n)\right) - \mu'_t(\mathscr{L}f_n)}. 
    \end{align*}
   Combining this with Proposition \ref{proposition:local-infinitesimal-rotation-property} implies that for all $n \in \N$
    \begin{align*}
        \limsup_{\varepsilon \downarrow 0}\abs{\tfrac{1}{\varepsilon}\left(\mu'_{t+\varepsilon}(f)-\mu'_t(f)\right) - \mu'_t(\mathscr{L}f)} 
        \leq 
        C \vertiii{f_n -f} + \norm{\mathscr{L}f-\mathscr{L}f_n}_\infty. 
    \end{align*}
    Now, we use that both terms on the right-hand side vanish as $n$ tends to infinity to conclude that
    \begin{align*}
        \frac{d}{d\varepsilon}\lvert_{\varepsilon=0}\mu'_{t+\varepsilon}(f) = \mu'_t(\mathscr{L}f). 
    \end{align*}
    Therefore, the infinitesimal rotation property can be extended from all local functions to all functions in $D(\Omega)$. 
\end{proof}

\subsection{The forward-backward construction}

\begin{proof}[Proof of Proposition \ref{proposition:time-derivative-forward-backward}]
    For $h \in \R$, with $\abs{h}$ sufficiently small, we have 
    \begin{align*}
        \tfrac{1}{h}&\left(\mu'_{s+h}(S(t-s-h)f)-\mu'_s(S(t-s)f)\right)
        \\\
        &=
        \tfrac{1}{h}\left(\mu'_{s+h}(S(t-s-h)f)-\mu'_s(S(t-s-h)f)\right)
        \\\
        &\qquad+
        \tfrac{1}{h}\left(\mu'_s(S(t-s-h)f)-\mu'_s(S(t-s)f)\right). 
    \end{align*}
    For the second summand, we can simply use the definition of the generator and dominated convergence to see that 
    \begin{align*}
        \lim_{h \to 0}\tfrac{1}{h}\left(\mu'_s(S(t-s-h)f)-\mu'_s(S(t-s)f)\right) = - \mu'_s(\mathscr{L}S(t-s)f). 
    \end{align*}
    Now, for the first summand, we use that $S(r)f \in D(\Omega)$ for all $r\geq 0$ and the infinitesimal rotation property to see that 
    \begin{align*}
        \lim_{h\to 0}\tfrac{1}{h}\left(\mu'_{s+h}(S(t-s)f) - \mu'_s(S(t-s)f)\right) = \mu'_s(\mathscr{L}S(t-s)f). 
    \end{align*}
    Moreover, for any function $g \in D(\Omega)$ we can apply Lemma \ref{lemma:lipschitz-continuity-of-rotation} to get 
    \begin{align*}
        \tfrac{1}{h}\abs{\mu'_{s+h}(g)-\mu'_s(g)} \leq C \vertiii{g}.
    \end{align*} 
    Hence, applying this estimate with $g=S(t-s-h)f-S(t-s)f$ and use the strong continuity of the semigroup plus the upgrading lemma to conclude that 
    \begin{align*}
        \tfrac{1}{h}&\abs{\mu'_{s+h}(S(t-s-h)f-S(t-s)f)-\mu'_s(S(t-s-h)f-S(t-s)f)} 
        \\\
        &\leq 
        C  \vertiii{S(t-s-h)f-S(t-s)f} \to 0
    \end{align*}
    as $h \to 0$. Putting all of these pieces together yields 
    \begin{align*}
        \lim_{h \to 0}\tfrac{1}{h}&\left(\mu'_{s+h}(S(t-s-h)f)-\mu'_s(S(t-s)f)\right)
        \\\
        &= 
        \mu'_s(\mathscr{L}S(t-s)f) - \mu'_s(\mathscr{L}S(t-s)f) 
        = 
        0,
    \end{align*}
    as we have claimed. 
\end{proof}

\subsection{Uniqueness of the equilibrium measure}\label{sec:proof-uniqueness-equilibrium-measure}

For the proof we need the following two results from the literature as technical helpers.  The first one tells us something about the structure of the set of translation-invariant extremal Gibbs measures for the specification $\gamma$. 

\begin{proposition}\label{prop_extr}
    For the extremal translation-invariant Gibbs measures for the specification $\gamma$ we have 
    \begin{align*}
        \text{ex}\mathscr{G}_\theta(\gamma) 
        =
        \left\{
            \mu_\varphi = \lim_{\Lambda \uparrow \Z^d}\gamma_\Lambda(\cdot \lvert \eta_\varphi) : \ \varphi \in \mathbb{S}^1
        \right\},
    \end{align*}
    where $\eta_\varphi$ is the configuration that is equal to $\varphi$ everyhwere. 
    Therefore, the elements of $\emph{ex}\mathscr{G}_\theta(\gamma)$ can be  uniquely labelled by the angle $\varphi \in [0,2\pi)$. 
\end{proposition}

For the classical XY model in $d\geq 3$ a structural result of this type was first shown in \cite[Section 3]{frohlich_spin_1983} and later extended to more general systems – including long-range models in dimension one and two – in \cite{pfister_ergodic_1987}. \medskip 

Additionally, recall that by Proposition~\ref{proposition:kernel-correspondence-discrete-continuous}, the discretisation actually induces a bijection between the translation-invariant extremal Gibbs measures of $\gamma$ and the translation-invariant extremal Gibbs measures of $\gamma'$, so we can also label the extremal Gibbs measures for the discrete specification. 

Now we can show that the equilibrium measure for the dynamics generated by $\mathscr{L}$ is indeed unique.  

\begin{proof}[Proof of Proposition \ref{proposition:uniqueness-equilibrium-measure}]
    First note that, by Proposition~\ref{proposition:full-rotation-property} the measure $\mu^*$ defined by 
    \begin{align*}
        \mu^* = \frac{1}{2\pi}\int_0^{2\pi}\mu_t' dt. 
    \end{align*}
    is indeed time-stationary for the dynamics, because it is a Cesàro average along a measure-valued trajectory, see for example \cite[Proposition I.1.8(e)]{liggett_interacting_2005}. 

    Now let $\nu' \in \calM_1(\Omega)$ be translation-invariant and time-stationary w.r.t.~the dynamics generated by $\mathscr{L}$. 
    Then, by  Proposition~\ref{proposition:attractor-property}, we have $\nu' \in \mathscr{G}(\gamma')$. Hence, by the extremal decomposition for Gibbs measures, see~\cite[Theorem 7.26]{georgii_gibbs_2011}, we know that $\nu'$ can be written as
    \begin{align*}
        \nu' = \int_{\text{ex}\mathscr{G}_\theta(\gamma')}\overline{\mu}w_{\nu'}(d\overline{\mu}),
    \end{align*}
    where $w_{\nu'}$ is a probability measure on the set $\text{ex}\mathscr{G}_\theta(\gamma')$. By combining Proposition~\ref{prop_extr} and Proposition~\ref{proposition:kernel-correspondence-discrete-continuous} we know that there is a measurable bijection 
    \begin{align*}
        b: \text{ex}\mathscr{G}_\theta(\gamma') \to [0,2\pi), \quad \mu' \mapsto \text{arg}(\mu'(\omega'_0)). 
    \end{align*}
    Hence, we can consider the image measure $v_{\nu'}$ of $w_{\nu'}$ under this mapping. Note that $v_{\nu'}$ is now a probability measure on $[0,2\pi)$. 
    The main idea is now to show that $v_{\nu'}$ is a translation invariant measure on $[0,2\pi)$ (with periodic boundary on $[0,2\pi)$ of course), because by the standard characterisation of the Lebesgue measure this directly implies $v_{\mu'} (dt) = \frac{1}{2\pi}\lambda(dt)$ and hence $\nu' = \mu^*$. 
    For the slightly technical details see \cite[Proposition 5.1]{jahnel_class_2014}.
\end{proof}

\section{Outlook and open problems}\label{section:outlook} 
Let us briefly mention some interesting related open problems and some ideas for future research. 
\subsubsection{Time-periodic behaviour in reversible systems}
It seems to be folklore that time-periodic behaviour can only happen in driven, i.e., irreversible, systems. But in the literature there is actually neither a robust heuristic reason, nor a rigorous mathematical argument that justifies this assumption. As a thought experiment, we would like to invite the reader to consider the following construction which could be interpreted as the dynamical counterpart to Dobrushin's construction of non-translation-invariant Gibbs measures in dimensions $d\geq3$. \medskip

As local state space, consider a (sufficiently fine) discretisation of the unit sphere, e.g., 
\begin{align*}
    \Omega_0 = \left\{\exp(i\varphi): \ \varphi=2\pi k/N, k=1,\dots,N-1 \right\} 
\end{align*}
and assume that the interaction is given by the (formal) Hamiltonian 
\begin{align*}
    \calH(\varphi) = -\sum_{x\sim y}\beta \cos(\varphi_x - \varphi_y). 
\end{align*}
For the transition rates of the system on all of $\Z^d$ consider Glauber-type rates which only allow jumps by one discrete step 
\begin{align*}
    c(x,\varphi, \pm) = \exp\left(\beta\sum_{x \sim y}\left[\cos\left(\varphi_x - \varphi_y \pm \frac{2\pi}{N}\right) - \cos(\varphi_x - \varphi_y)\right] \right), 
\end{align*}
where $c(x, \varphi, +)$ is the rate of performing a clockwise rotation by an angle of $2\pi/N$ at site $x$ and $c(x,\varphi,-)$ is the rate of performing a counter-clockwise rotation. Clearly, the Gibbs measures w.r.t.~the Hamiltonian $\calH(\cdot)$ are reversible measures for the process with these rates as one can see by checking the detailed-balance equations. \medskip 

In finite volumes $\Lambda \Subset \Z^d$ we will consider a time-dependent (but periodic) generator which is constructed as follows. 
For $\omega \in [0,2\pi)$ denote by $\overline{\omega}$ the configuration in $(\mathbb{S}^1)^{\Z^d}$ that is constantly equal to $\exp(i\omega)$ and consider the time-periodic map 
\begin{align*}
    [0, \infty) \ni t = \omega(t) + 2\pi \lfloor t/(2\pi)\rfloor \mapsto \overline{\omega}(t) := ( t \mod 2\pi ) \in (\mathbb{S}^1)^{\Z^d}.
\end{align*}
Then we can define the time-dependent rates of our finite-volume processes by 
\begin{align*}
    c^\Lambda_t(x, \eta, \pm) 
    = 
    \frac{\exp\left(\beta \calH(\eta^{x,\pm}_{\Lambda} \overline{\omega}(t)_{\Lambda^c}\right)}{\exp(\beta \calH(\eta_{\Lambda}\overline{\omega}(t)_{\Lambda^c})},
\end{align*}
where $\eta^{x,\pm}$ is the configuration which agrees with $\eta$ everywhere except at site $x\in \Lambda$ where it is equal to $\eta_x \pm 2\pi/N$ (so we again rotate the spin by one discrete step).
This gives us a (time-dependent) generator of the form 
\begin{align*}
    \mathscr{L}_t^\Lambda f(\eta)= \sum_{x \in \Lambda}c^\Lambda_t(x, \eta, \pm)[f(\eta^{x,\pm})-f(\eta)], 
\end{align*}
which gives rise to a flow $(S^\Lambda_{t,s})_{0 \leq s \leq t}$ on the space of continuous functions $C(\Omega)$. \medskip

As $\Lambda \uparrow \Z^d$, the rates of the finite-vlume dynamics converge uniformly to the rates of the reversible dynamics generated by $\mathscr{L}$, so one might expect to have convergence of the flows $(S^\Lambda_{s,t})_{0\leq s \leq t}$ to the semigroup $(S(s))_{s \geq 0}$ via a time-dependent version of the Trotter–Kurtz theorem, see e.g., \cite[Proposition I.2.12]{liggett_interacting_2005}. 

The most interesting part comes now. In finite volumes $\Lambda$, the time-reversal symmetry of the rates is broken by the rotating boundary conditions and one is essentially in the situation of classical theory for linear ODEs with time-periodic coefficients with periodic $2\pi$.  One can show that there exists a unique periodic solution with period $2\pi$ and that this solution is moreover attractive, see e.g., \cite[Proposition 7.3]{bertini_level_2018}. \medskip

The question is now whether this time-periodic orbit survives the thermodynamic limit $\Lambda \uparrow \Z^d$, i.e., if it converges to a non-trivial time-periodic orbit for the reversible dynamics generated by $\mathscr{L}$. This would give us a breaking of the time-reversal symmetry for sufficiently strong interactions, in which the system – when started in a specific initial condition – never forgets the direction of time, despite its reversible nature. \medskip

Note, that in \cite{jahnel_long-time_2023} it was shown that in dimensions $d=1,2$ and under short-range and strong irreducibility assumptions, reversible interacting particle systems cannot exhibit time-periodic behaviour. So if one wants to work in dimensions $d=1,2$, then one definitely needs to violate one of these assumptions to make the construction work. In higher dimensions, there are – at least for the moment – no such theoretical obstructions. 
\subsubsection{Extensions of the results of Mountford and Ramirez–Varadhan}
Another possible direction for future research is to try to extend the results from \cite{mountford_coupling_1995} and \cite{ramirez_relative_1996} from finite-range systems to systems with unbounded but reasonably fast decaying interactions in dimension $d=1$ or in general to two-dimensional systems. Since both proof strategies use the one-dimensional and finite-range structure quite explicitly, one needs to come up with some novel ideas to overcome some challenging technical problems. However, one could learn a lot about the approach to equilibrium of interacting particle systems, so it seems like a hard but worthwhile endeavour. 

\subsubsection{Time-periodic behaviour for systems with explicit finite-range interactions}
One big drawback of the construction in this article and in \cite{jahnel_class_2014} is of course that the rates are only semi-explicit and certainly not finite range. To gain a better understanding of the emergence of time-periodic behaviour, one would like to construct an, as simple as possible, non-degenerate interacting particle system that exhibits time-periodic behaviour and there have been some proposals for such systems, see \cite[Chapter 1.8]{swart_course_2022}. However, there has been little success in actually verifying the existence of time-periodic orbits. This is mostly due to the fact that one lacks a precise description of how these orbits look, in other words, one does not have a clear description of the measures on them. 

\subsection*{Acknowledgements}
The authors would like to thank the anonymous referee for their insightful feedback that helped to substantially improve this manuscript.
Additionally, the authors would like to thank Aernout van Enter, Christof Külske, and Jan Swart for helpful discussions and comments during the preparation of this manuscript.  
The authors acknowledge the financial support of the Leibniz Association within the Leibniz Junior Research Group on \textit{Probabilistic Methods for Dynamic Communication Networks} as part of the Leibniz Competition.

\subsection*{Data availability} Data sharing not applicable to this article as no data sets where generated or analysed during the current study. 

\subsection*{Statements and Declarations}
The authors have no relevant financial or non-financial interests to disclose.

\bibliography{references}
\bibliographystyle{alpha}
\end{document}